\documentstyle[epsfig, 12pt]{article}

\newcommand{\qed}{{$\Box$}}

\newcommand{\ba}{\begin{array}}
\newcommand{\ea}{\end{array}}

\newtheorem{lemma}{Lemma}
\newtheorem{hypo}{Assumption}
\newcommand{\bhypo}{\begin{hypo}}
\newcommand{\ehypo}{\end{hypo}}
\newtheorem{defi}{Definition}
\newcommand{\ble}{\begin{lemma}}
\newcommand{\ele}{\end{lemma}}
\newcommand{\bde}{\begin{defi}}
\newcommand{\ede}{\end{defi}}

\newtheorem{prop}{Proposition}
\newcommand{\epr}{\end{prop}}
\newcommand{\bpr}{\begin{prop}}
\newtheorem{teo}{Theorem}
\newcommand{\bth}{\begin{teo}}
\newcommand{\eth}{\end{teo}}
\newtheorem{rema}{Remark}
\newcommand{\bre}{\begin{rema}}
\newcommand{\ere}{\end{rema}}
\newcommand{\ee}{\end{equation}}
\newcommand{\be}{\begin{equation}}

\def\sq#1{\overline{\underline{\left|\matrix{\cr~~ #1~~\cr\cr}\right|}}}

\begin{document}

\leftline{Running head: {\bf Solving the Sixth Painlev\'e Equation}}
\vskip 0.5 cm

\centerline{\Large\bf Solving the Sixth Painlev\'e Equation:}
\vskip 0.3 cm 
\centerline{\Large\bf  Towards the Classification of all}
\vskip 0.3 cm 
\centerline{\Large\bf  
   the Critical Behaviors and the Connection Formulae  }

\vskip 0.5 cm
\centerline{{\Large Davide Guzzetti}}
~~~~~~~~~~~~~~~~~~~~~~~~~~~~~~~~~~~~~~~~~~~~~~~~~~~~~~~~~~~~~~~~~~~~\footnote{
Korea Institute of
    Advanced Study KIAS, Hoegiro 87(207-43 Cheongnyangni-dong), Dongdaemun-gu, Seoul 130-722, Korea. E-mail: guzzetti@kias.re.kr~~~ Tel: +82-2-958-3861 ~~~Fax: +82-2-958-3786.

Also: International School of Advanced Studies SISSA/ISAS, Via Bonomea 265, 34136 Trieste, Italy.
}


\vskip 0.5 cm 
\noindent
{\bf Abstract:} The critical behavior of a three real parameter class of solutions of the sixth Painlev\'e equation is computed, and parametrized in terms of monodromy data of the associated $2\times 2$ matrix linear Fuchsian system of ODE. The class may contain solutions with poles accumulating at the critical point.   The study of this class  closes a gap in the description of the transcendents in one to one correspondence with the monodromy data. These transcendents  are reviewed in the paper.  Some formulas that relate the monodromy data to the critical behaviors of the  four real (two complex) parameter class of solutions are missing in the literature, so they are computed here.     A computational procedure to write
 the full expansion of the four and three real parameter class of solutions is proposed.

\section{Introduction}

 The history, importance and applications of the Painlev\'e equations  have been
widely discussed in the literature and assumed to be known
 (for a review, see \cite{Its}). The equation PVI is: 
 $$
{d^2y \over dx^2}={1\over 2}\left[ 
{1\over y}+{1\over y-1}+{1\over y-x}
\right]
           \left({dy\over dx}\right)^2
-\left[
{1\over x}+{1\over x-1}+{1\over y-x}
\right]{dy \over dx}
$$
$$
+
{y(y-1)(y-x)\over x^2 (x-1)^2}
\left[
\alpha+\beta {x\over y^2} + \gamma {x-1\over (y-1)^2} +\delta
{x(x-1)\over (y-x)^2}
\right]
,~~~~~\hbox{(PVI)}.
$$ 
The general solution has no movable essential singularities or branch points, which are possibly located only at the 
 {\it critical points} $x=0,1,\infty$. The
behavior of a solution when $x\to 0,1,\infty$, is called {\it critical
  behavior}. The other movable singularities are poles.  
The absence of movable 
critical points means that a solution
can be meromorphically extended to the
universal covering of a punctured complex sphere, determined  only by
the equation. Thus PVI shares a fundamental property of the linear
equations defining classical transcendental functions.

Following the review \cite{Its}, an expression is called  {\it
  explicit} when it is given in terms of a finite
  algebraic 
  combination of elementary and elliptic functions, and a finite
  number of contour integrals (and quadratures) of these
  functions. Classical linear special functions admit explicit
  representations.  The general solution of a Painlev\'e equation  does not, as it is
  proved by H.Umemura in  \cite{Umemura}. Therefore, it is not a classical function. It is a new function,  called a {\it Painlev\'e transcendent}.\footnote{
 H. Umemura proved the of irreducibility of the Painlev\'e equations \cite{Umemura} \cite{Umemura1} \cite{Umemura2}. The term "explicit" expression is equivalent to the notion of classical function.  
Following \cite{Umemura}, a function is called {\it classical} if it  is given in terms of a finite iteration of {\it permissible operations} applied to rational functions.  They are the derivation,  rational combination  (sum, product, quotient), algebraic combinations (the expression is a root of a polynomial whose coefficients are rational functions (and then, after iteration,  classical functions)), contour integrals and quadratures, solution of a linear homogeneous differential equation whose coefficients are rational functions (or classical functions, after iteration),        a solution of an algebraic differential equation of the first order whose coefficients are rational functions (or classical functions), 
composition with abelian 
functions (the expression is  $\varphi(f_1(x),...,f_n(x))$, where $f_1$,...,$f_n$ are rational or classical functions, and $\varphi:{\bf C}^n/ \Gamma \to {\bf C}$ is meromorphic, $\Gamma$ is a lattice). The reader may note that the elementary transcendental functions are classical functions (they are the algebraic functions, or a  function which is obtained from an algebraic function by integration (like the exponential, the trigonometric and hyperbolic functions), or the inverse of such an integral (like the logarithm, the elliptic functions, etc)).
Umemura proved in \cite{Umemura} that the general solution of a Painlev\'e equation is not a classical function.      
 H.Watanabe \cite{watanabe} applied the argument to PVI, and showed that a solution of PVI is either algebraic, or
solves a Riccati equation (one-parameter family of classical solutions), or it is not a classical  function.  All the algebraic solutions were classified in \cite{DM} when
  $\beta=\gamma=0$,
 $\delta={1\over 2}$, and then in \cite{Lisovyy} for the general PVI. 
}

\vskip 0.2 cm 
 Following  \cite{Its} (page 8), 
 {\it solving} PVI means:  
~i) Determine 
the {\it explicit} critical behavior of the transcendents 
at the  critical points. Such a behavior must be given by  an
 explicit formula in terms of 
two integration constants. ~ii) Solve the {\it connection
  problem}, namely: find the {\it explicit} 
relations among couples of integration
constants  at different critical points. The above i) and ii) are the
problem of {\it global analysis} of the equation. 
 Solution of  i) and ii) means that a Painlev\'e transcendents can be
efficiently used in applications 
as it is the case for special functions.   
It was thought that the global analysis is possible only for linear equations, namely only  for classical linear special functions. But the {\it method of monodromy preserving deformations}  has made the global analysis possible also for Painlev\'e equations.\footnote{A more restrictive definition of ``solving'' should include 
the distribution of the poles (movable singularities) of the
 transcendents. This problem for PVI is still open (in \cite{D4}, 
 the behavior on the universal covering of a critical point is analyzed and it is shown that if the poles exist, they are distributed in spirals converging to the critical point).}

\vskip 0.2 cm

 The critical behaviors  for a two complex (four real) parameter class of solutions  were  computed  and parametrized in terms of monodromy data of an associated Fuchsian system of ODEs, by Jimbo  in \cite{Jimbo}. Jimbo's paper is the foundation of all the works on PVI  based on the method of monodromy preserving deformations which have followed. 
  Some authors have  determined critical
 behaviors not included in Jimbo's  class, with different methods. Among them, the works of S.Shimomura (the results are summarized
 in \cite{IKSY}) and   A.D.Bruno, I.V. Goryuchkina (\cite{Bruno}
 \cite{Bruno1} \cite{Bruno2} \cite{Bruno3} \cite{Bruno6}) are local approaches, which do not determine the connection
 formulae, but essentially
 determine all the critical behaviors or asymptotic expansions. In \cite{D4} \cite{D3} \cite{D2} \cite{D1}, in the framework of the method of monodromy preserving deformations,   critical behaviors not included in Jimbo's  class are constructed and parametrized in terms of associated monodromy data.    Accordingly, the transcendents can be cassified into a few  classes (one beeing that of Jimbo's), depending on  their local behavior and their correspondence with  sub spaces of the space of the associated monodromy data. This fact is reviewed below, in subsection \ref{review}. The last class (a three real parameter cass of solutions)  has not been  studied yet, and it is studied here. Its critical behavior and parametrization in terms of monodromy data is given in the present paper. Together with it, in the paper a complete review of the parametrization of critical behaviors in terms of monodromy data is given for all the  four and three real parameter solutions. Some formulas missing in the literature are computed.   
 
According to the above definition of ``solving'', 
 the paper by Jimbo \cite{Jimbo} and \cite{D4} \cite{D3} \cite{D2} \cite{D1}, together with the present paper,  "solve"
 PVI. Namely,   the critical behaviors and the parametrization in terms of monodromy data have been found for of all the  transcendents that are in  one-to-one correspondence with points in the space of the associated monodromy data. 

Before stating the results of the paper, we give a review of the critical behaviors at $x=0$.


\subsection{A Review of Critical Behaviors}
\label{review}
\vskip 0.2 cm 
In the following, let $|\arg x|<\pi$, $|\arg (1-x)|<\pi$, so that all functions of $x$ will be understood as  $x$-branches.  
According to \cite{JMU},  
 PVI is the condition of isomonodromy deformation for a $2\times 2$
 fuchsian system with four singularities $0,x,1,\infty$:  
\be
   {d\Psi\over d\lambda}=A(x,\lambda)~\Psi,~~~~~
A(x,\lambda):=\left[ {A_0(x)\over \lambda}+{A_x(x) \over \lambda-x}+{A_1(x)
\over
\lambda-1}\right],~~~\lambda\in{\bf C}.
\label{SYSTEM}
\ee
The traces of the matrices are zero, and the eigenvalues are fixed by PVI. This facts are reviewed in  Section \ref{REVIEW}. A fundamental solution $\Psi$ has  branch points in $\lambda = 0,x,1$. Fix a base point and a base of loops $\Gamma$  like in figure \ref{figure1}.   When $\lambda$ goes around a small loop around a branch point, the  fundamental solution transforms like $\Psi \mapsto \Psi M_i$, $i=0,x,1$.  The $2\times 2$ matrices    $M_0$, $M_x$, $M_1$ are called the monodromy matrices  of the fundamental solution. 

Given PVI (namely, given $\alpha$, $\beta$, $\gamma$ and $\delta$), there is a one-to-one
 correspondence between a triple of monodromy martices, associated to the base of loops $\Gamma$,  and a branch of a PVI
 transcendent. This happens  in the generic case (which will be made  precise in Section \ref{REVIEW}). A branch is uniquely identified by the monodromy data associated to the  basis of loops $\Gamma$: 
$$
y(x)=y(x; \hbox{Tr$M_0$,    Tr$M_x$, Tr$M_1$, Tr$M_0M_x$, Tr$M_xM_1$, Tr$M_0M_1$})
$$
 This will be precisely explained  in Section \ref{REVIEW}. Here it is enough to understand that the critical behavior at a critical point depends on two integration constants (which in general are 4 real parameters, but in sub cases they may  reduce to 3 or 2 real parameters). The parametrization of the integration constants in terms of monodromy data  uniquely identifies the branch of the transcendent.  The expicit parametrization will be given in Section \ref{conection}, for the 4-real parameter and 3-real parameter branches. 

\vskip 0.2 cm 

 What kind of  critical behaviors we may expect at $x=0,1,\infty$ depends 
on the values  of Tr$M_0M_x$, Tr$M_xM_1$, Tr$M_0M_1$ respectively. For example, the  type of behavior at $x=0$  (for example, a two real parameter solution with logarithmic behavior, or a solution of Jimbo's, etc)
is decided by the value of Tr$(M_0M_x)$. 

Here the classes of critical behaviors are reviewed, corresponding to monodromy groups which have the property of being  irreducible, and such that they are in one to one correspondence with branches of PVI transcendents (namely, none of the monodromy matrices $M_0,M_x,M_1,M_1M_0M_x$ is the identity).

\vskip 0.2 cm 

$\diamond$) [{\bf Small power type behaviors (Jimbo). 4-real parameters}:]  
 M.Jimbo was the first to determine the critical behaviors for a wide
 class of transcendents. In \cite{Jimbo} he proved that PVI admits 
 solutions with branches behaving as follows: 
\be
\label{INTROO1}
y(x)=
\left\{\matrix{
a_{0}x^{1-\sigma_{0}}(1+\delta_0(x)),~~~x\to 0
\cr\cr
1-a_1(1-x)^{1-\sigma_{1}}(1+\delta_1(1-x)),~~~x\to 1
\cr\cr
a_\infty x^{\sigma_{\infty}}(1+\delta_\infty(x^{-1})),~~~x\to \infty
}
\right.
\ee
where $a_{i},\sigma_i\in{\bf C}$ are integration constants such that:
$$
a_i\neq 0,~~~~~0<\Re\sigma_i<1. 
$$
$\delta_i(\zeta)$ are higher order terms, $\delta(\zeta)=
O(\hbox{max}\{|\zeta|^{\Re\sigma},|\zeta|^{1-\Re \sigma}\})$. 
  Jimbo determined the parametrization of the couples
  $(a_0,\sigma_0)$ $(a_1,\sigma_1)$, $(a_\infty,\sigma_\infty)$ in terms of monodromy data.  {\it The parametrization identifies the specific branch}. 
In particular he proved that: 
\be
  2\cos (\pi \sigma_0) = \hbox{Tr}(M_0M_x),~~~ 2\cos (\pi \sigma_1) =
  \hbox{Tr}(M_xM_1),~~~  2\cos (\pi \sigma_\infty) = \hbox{Tr}(M_0M_1)
\label{12novembre2010}
\ee
The restriction on $\Re \sigma_i$ means that the solutions correspond to the following subspace of the space of monodromy matrices: 
$$
  \hbox{Tr}(M_iM_j)\not \in (-\infty,-2]\cup [2,\infty). 
$$

  For special values of $\sigma$ the above behaviors are modified. For example,  for $x\to 0$, one has (see \cite{D2}, plus section \ref{FULL} and section \ref{maldipiedi} of the present paper): 
$$
y(x)={\sqrt{-2\beta}\over \sqrt{-2\beta}+\sqrt{1-2\delta}}x\mp {r\over  \sqrt{-2\beta}+\sqrt{1-2\delta}} ~x^{1+\sigma} +O(x^2),~~~~~\sigma=\pm(\sqrt{-2\beta}+\sqrt{1-2\delta}),
$$
$$
y(x)={\sqrt{-2\beta}\over \sqrt{-2\beta}-\sqrt{1-2\delta}}x\mp {r\over  \sqrt{-2\beta}-\sqrt{1-2\delta}} ~x^{1+\sigma}+O(x^2) ,~~~~~\sigma=\pm(\sqrt{-2\beta}-\sqrt{1-2\delta}),
$$
Here $r\in {\bf C}$ is the integration constant and the condition $-1<\Re\sigma<1$ must hold.

\vskip 0.2 cm 
Not only in Jimbo's case, but in general,  the  critical behavior of 
  $y(x)$ is  decided  by three constants  
   $\sigma_0$, $\sigma_1$, $\sigma_\infty$, 
  determined 
   by (\ref{12novembre2010})  plus the 
   conditions $0\leq \Re \sigma_i \leq 1$. 

Below, 
behaviors are given only for  $x\to 0$ (arg$(x)$ bounded). We  denote 
$\sigma:=\sigma_0$. 
The other critical points $x=1,\infty$ will be
described in the paper.

\vskip 0.2 cm

$\diamond$) [{\bf Sine-type oscillatory behaviors. 3 real parameters}:]  If $\Re\sigma=0$,
  the critical behavior follows from Jimbo's results (see
  Appendix I) and the equivalent method of \cite{D2}. There exist a transcendent with a branch at $x=0$ behaving as follows:  
 \be
\label{INTROO2}
y(x)= x\left\{
iA\sin\left(
i\sigma\ln x +\phi
\right)~+B+\delta^{*}(x)
\right\},~~~\delta^{*}(x)=O(x),~~~x\to 0
\ee
$$
\sigma,~\phi \hbox{ integration constants. }~~
B ={\sigma^2-2\beta-1+2\delta \over 2\sigma^2},~~~
A^2+B^2=-{2\beta\over \sigma^2}.  
$$
 In this case:  $$
2\cos\pi\sigma=\hbox{Tr}(M_0M_x)>2.
$$
The parametrization of $\sigma$ and $\phi$ in terms of monodromy data uniquely identifies the branch.

\vskip 0.2 cm 
$\diamond$) [{\bf Log-type behaviors. 2 real parametes}:]  If $\sigma=0,1$, namely: 
$$
 \hbox{Tr}(M_0M_x)=\pm 2,
$$ 
 There are transcendents  with logarithmic branches (see \cite{Jimbo} formula (1.9$)^\prime$, and \cite{D2} \cite{D1}). In \cite{D2} \cite{D1}, the branches  are written as follows. 
When $\sigma=0$, Tr$(M_0M_x)=2$: 
$$
  y(x)= x\left[
{1+2\beta-2\delta\over 4} \left(
\ln x + {4r +2\sqrt{-2\beta} \over 2\delta -2\beta -1}
\right)^2 +{2\beta \over 2\beta +1-2\delta}
\right]+O(x^2\ln^3 x),~~~2\beta\neq 2\delta-1;
$$
$$
y(x)= x(r\pm \sqrt{-2\beta} \ln x) +O(x^2\ln^2 x),~~~2\beta= 2\delta-1.
$$
For the second solution, the subgroup $<M_0,M_x>$ is reducible. $r$ is the integration constant. 
\vskip 0.2 cm
\noindent
 When $\sigma=1$, Tr$(M_0M_x)=-2$: 
$$
 y(x)= {2\over (\gamma-\alpha)\ln^2 x} \left[
1+ {4r +\sqrt{8\alpha} \over \gamma -\alpha} {1\over \ln x} +O\left(
{1\over \ln^2 x}
\right)
\right],~~~\alpha\neq \gamma;
$$
$$
 y(x)= {1\over  \pm\sqrt{2\alpha}~\ln x}\left[1\mp{r\over \sqrt{2\alpha} ~\ln x} +O\left(
{1\over \ln^2 x}\right) \right],~~~\alpha = \gamma. 
$$
For the second solution, the subgroup $<M_0M_x,M_1>$ is reducible. 
$r$ is the integration constant.  Its parametrized in terms of monodromy data is in \cite{D1}. This identifies the branch.

\vskip 0.2 cm 
$\diamond$) [{\bf Taylor expansions. 2 real parameters}:]   Solutions with branches which admit a Taylor expansions at a critical point are studied in \cite{D2}, \cite{Kaneko}. According to \cite{D2}, such expansions (which are convergent for small $|x|$ by the argument of  \cite{Kaneko})  are the following 1), 2), 3) below: 

\vskip 0.3 cm 
\noindent
1) Degenerate solutions $y=0,x,1$. 

\vskip 0.3 cm 
\noindent
2) The  {\it Basic Expansions} i), ii), iii) below:. 

\vskip 0.2 cm 
i) When $\alpha\neq 0$ and $\sqrt{2\alpha}\pm\sqrt{2\gamma}$ is not integer: 
$$
 y(x)= {\sqrt{\alpha}\pm\sqrt{\gamma}\over \sqrt{\alpha}}\mp
{
\sqrt{\gamma}~\Bigl[(\sqrt{2\alpha}\pm\sqrt{2\gamma})^2-2\delta+2\beta\Bigr]
\over 
2\sqrt{\alpha}((\sqrt{2\alpha}\pm\sqrt{2\gamma})^2-1)
}
~x+ \sum_{n=2}^\infty c_n(\sqrt{\alpha},\pm\sqrt{\gamma},\beta,\delta)~x^n
$$
$$\hbox{Tr}(M_0M_x)=-2\cos \pi (\sqrt{2\alpha}\pm\sqrt{2\gamma}).$$

\vskip 0.2 cm 
ii) When $\alpha\neq 0$, but $\sqrt{2\alpha}\pm\sqrt{2\gamma}= 1$ and $1-2\delta+2\beta=0$:
$$
 y(x) = \pm{1\over \sqrt{2\alpha}} + rx +  \sum_{n=2}^\infty c_n(r;\sqrt{\alpha},\beta)~x^n,~~~~~r\in{\bf C}
$$
$$\hbox{Tr}(M_0M_x)=2.$$

\vskip 0.2 cm 
iii) When $\alpha= 0$ and $\sqrt{\alpha}\pm\sqrt{\gamma} =0$:
$$
 y(x)= r+(1-r) (\delta-\beta)x+\sum_{n=2}^\infty c_n(a;\beta,\delta)~x^n,~~~~~r\in{\bf C}
$$
$$\hbox{Tr}(M_0M_x)=-2.$$

\vskip 0.2 cm
\noindent
In all i), ii), iii) above, the subgroup $<M_0M_x,M_1>$ is reducible. 
The square roots $\sqrt{\alpha}$, $\sqrt{\gamma}$ have arbitrary sign. The coefficients are rational functions of their arguments. 
 The parametrization of $r$ in terms of monodromy data is in \cite{D2}. It uniquely identifies the branch. 
\vskip 0.3 cm 
\noindent
3) All the expansions obtained from 2) by the birational transformations  of PVI that do not change $x$. For example, the bitrational transformation (\ref{infi}) gives: 

\vskip 0.2 cm 
i) When $\beta\neq 0$ and $\sqrt{-2\beta}\pm\sqrt{1-2\delta}\neq 0$:

$$
y(x)= {\sqrt{-2\beta}~x\over \sqrt{-2\beta}\pm\sqrt{1-2\delta}}~\pm{
\sqrt{-2\beta}\sqrt{1-2\delta}\Bigl[(\sqrt{-2\beta}\pm\sqrt{1-2\delta})^2+2\gamma-2\alpha-1\Bigr]~x^2
\over 
2(\sqrt{-2\beta}\pm\sqrt{1-2\delta})^2\Bigl[(\sqrt{-2\beta}\pm\sqrt{1-2\delta})^2-1\Bigr]
}~
+
$$
$$
+\sum_{n=3}^\infty b_n(\alpha,\sqrt{\beta},\sqrt{1-2\delta},\gamma)x^n
$$
$$\hbox{Tr}(M_0M_x)=-2\cos\pi(\sqrt{-2\beta}\pm\sqrt{1-2\delta}).
$$

\vskip 0.2 cm 
ii) When $\beta\neq 0$ but $(\sqrt{-2\beta}\pm\sqrt{1-2\delta})^2=1$ and $\alpha=\gamma$: 
$$
y(x)=\pm \sqrt{-2\beta} ~x~+r~x^2~+\sum_{n=3}^\infty b_n(r;\sqrt{\alpha},\sqrt{\beta})x^n,~~~~~r\in{\bf C}
$$
$$\hbox{Tr}(M_0M_x)=-2.$$

\vskip 0.2 cm 
iii) When $\beta=1-2\delta=0$:
$$
y(x)=rx+{r(r-1)\over 2}(2\gamma-2\alpha-1)x^2+\sum_{n=3}^\infty b_n(r;\alpha,\gamma)x^n,~~~~~r\in{\bf C}
$$
$$\hbox{Tr}(M_0M_x)=2.$$

\vskip 0.2 cm 
\noindent
In all I), II), III), the subgroup $<M_0,M_x>$ is reducible. 
\

\vskip 0.3 cm 
$\diamond$)  [{\bf Inverse sine-type oscillatory behaviors. 3 real parameters}:] 
The above results ``solve'' PVI for all the values of Tr$(M_iM_j)$, except
for the case $
\hbox{Tr}(M_iM_j) <-2$, namely the case when $\Re\sigma_i=1$. 
This case is studied in the present paper. The result, at $x=0$, is Proposition \ref{prop3}: there exist transcendents with a branch at $x=0$ having the following behavior: 
\be
\label{INTROO3}
y(x)= {1\over 
iA\sin\Bigl(i(1-\sigma_0)\ln x +\phi_0\Bigr)+B+\delta^*_0(x)},~~~\delta^*_0(x)=O(x),~~~x\to 0
\ee
$$
\sigma_0,~\phi_0 \hbox{ integration constants. }~~
B =
 {\Im\sigma_0^2+2\gamma-2\alpha\over 2
	\Im\sigma_0^2},~~~
A^2+B^2= -{2\alpha\over (\Im\sigma_0)^2}
$$
In this case:
$$
  2\cos(\pi \sigma_0)=\hbox{Tr}(M_0M_x)<-2,~~~~~\Re\sigma_0=1.
$$
The parametrization of $\sigma$ and $\phi$ in terms of monodromy data uniquely identifies the branch.

\vskip 0.3 cm 
\noindent
 A similar classification holds at $x=1,\infty$. Note that a solution with a  behavior falling in one class at a critical point, may have a behavior of a different type  at another critical point,  depending on  the values of Tr$(M_iM_j)$.


\subsection{Results of the Paper}

\noindent
The relevant results of this paper are the following three  points.  

\vskip 0.3 cm 
{\bf 1)} In this paper  PVI is solved in the missing case $
\hbox{Tr}(M_jM_k) <-2$, $\Re\sigma_i=1
$. 
Precisely: 

\vskip 0.2 cm 
 -- The critical behaviors when $x\to 0, 1,\infty$, with arg$x$
 and arg$(1-x)$ bounded is computed. Let PVI be given, and let the monodromy data be given (such that the
 one-to-one correspondence holds true). Let $x\to 0$
 inside a sector.  Let  
 $2\cos\pi\sigma_0=\hbox{Tr}(M_0M_x)<-2$, 
$\Re\sigma_0=1$.  The solution
 corresponding to these monodromy data has  the critical behavior (\ref{INTROO3}) [{\bf Inverse sine-type oscillatory behaviors}, Proposition \ref{prop3}].  Observe that:  
$$
\sin\Bigl(i(1-\sigma_0)\ln x +\phi\Bigr)= \sin\Bigl(\Im\sigma_0 \ln x +\phi\Bigr)
$$
gives a purely oscillating contribution when $x\to 0_+$. 
The above behavior also predicts the occurrence of poles close to
$x=0$, 
when the
denominator vanishes. This is the reason why the
correction $\delta^*_0(x)$ in the denominator must be kept. Namely, one  cannot write 
$y(x)= \{
iA\sin\Bigl(i(1-\sigma_0)\ln x +\phi_0\Bigr)+B\}^{-1}(1+O(x))$, because
this would affect the position of the poles. An example which makes
this point clear is the Picard-type solution (see Appendix II,
solution (\ref{PippoC})): 
$$
y(x)= {1+O(x)\over
\sin^2\left( {\Im\sigma\over 2} \ln x + \phi + {\Im\sigma\over 2} 
{F_1(x)\over F(x)} \right)} +O(x),  ~~~x\to 0
$$
where $F(x)$ and  $F_1(x)$ are the hypergeometric-like functions
(\ref{FffF}) and (\ref{FffF1}). 
The poles close to $x=0$ are determined by the solutions of 
 ${\Im\sigma\over 2} \ln x + \phi + {\Im\sigma\over 2} 
{F_1(x)\over F(x)}=k\pi$, $k\in {\bf Z}$, which
lie in a neighborhood of $x=0$. The distribution of poles in the
general case will be
studied in another paper.

\vskip 0.2 cm
 --  Connection Problem. The parametrization of the critical behavior in terms of monodromy data is given. In section \ref{conection}, Proposition
 \ref{prop6}  the critical behaviors at the 
three critical points $x=0,1,\infty$ is given, when $\Re
 \sigma_i=1$, $i=0,1,\infty$.  They are as follows: 
$$
y(x)= 1-{1\over 
iA_1\sin\Bigl(i(1-\sigma_1)\ln (1-x)
+\phi_1\Bigr)+B_1+\delta^*_1},~~~\delta^*_1=O(1-x),~~~x\to 1
$$
$$
y(x)= {x\over 
iA_\infty\sin\Bigl(i(\sigma_\infty-1)\ln x
+\phi_\infty\Bigr)+B_\infty+\delta^*_\infty(x)},~~~\delta^*_\infty=O\left({1\over
  x}\right),~~~x\to \infty
$$
The coefficients $A_1,A_\infty,B_1,B_\infty$ are given  in terms of
 $\alpha,\beta,\gamma,\delta$  in Proposition \ref{prop6}. 
 In Proposition \ref{prop6}  the integration constants $\phi_0,\phi_1,\phi_\infty$ are also given 
as functions of the coefficients
 of PVI and of the monodromy data Tr$(M_jM_k)$ (see (\ref{calR})). This parametrization fixes the  branches. 
  Conversely, in Proposition \ref{prop7},  the  formulae  which
express Tr$(M_jM_k)$ as functions of the coefficients of PVI and of
the integration constants are given. See formulae (\ref{COPPIAA1}).   
 In this way, one is able to compute  any of 
 the couples $(\sigma_0,\phi_0)$,  $(\sigma_1,\phi_1)$,
 $(\sigma_\infty,\phi_\infty)$ as a function of another. 
This {\it solves the connection problem}.

\vskip 0.2 cm 
 The author already studied the case 
Tr$(M_iM_j)<-2$ in \cite{D3}, \cite{D4}, with the elliptic
 representation. But the critical behavior obtained was  
 $y=\{\sin^2( {\Im \sigma\over 2} \ln x +\psi(x))+O(x)\}^{-1}$,
 where $\psi(x)= \sum_{n\geq 0}\psi_n x^{-in\Im\sigma }$
 is an oscillatory function.  The same behavior
 follows
 from the results of  Shimomura (\cite{IKSY}, chapter 4, section
 2). Unfortunately, the  function
 $\psi(x)$ in the sine 
 makes the formula uncomputable. The meaning of the result of the
 present paper  is 
 that {\it  $\psi(x)$ has been brought out } of the $\sin(~..~)$ and
 computed. 
The behavior of \cite{D3}, \cite{D4}, \cite{IKSY} of course
 must coincide with that of the present
 paper. This is possible because  one can always  write 
$iA\sin(\nu\ln x +\phi)+B$ (where $\nu\in{\bf
 R}$, $\phi\in{\bf C}$) as  $\sin^2({\nu
\over 2}\ln
 x +\sum_{n\geq 0} \psi(x))$, where $\psi(x)$ is an oscillating
 function (not vanishing for $x\to 0$) 
 computable in an elementary way. If $\psi(x)$ can
 be expanded in series in a suitable domain, then the series turns out
 to be necessarily of the form $\psi(x)= \sum_{n\geq 0}\psi_n x^{-in\nu }$,
 $\psi_n\in{\bf C}$. See Appendix II, subsection \ref{oscilla}, for the details.

\vskip 0.2 cm 
 It is to be cited the paper   \cite{Bruno6}, where all the   
asymptotic expansions are obtained with a power
 geometric technique \cite{Bruno5}. This technique  does not allow to
 solve the connection problem. In  \cite{Bruno6}, formula (7), one finds 
 an  expansion that, in the notation of the present paper, becomes   $y(x)=
 1/[iA\sin(i(1-\sigma)\ln x +\phi)+B] +\sum_{\Re s\geq 1} c_s x^s $. 
The absence of a term $\delta^*(x)=O(x)$ in the
 denominator, which is essential to determine the position of the
 poles, means that the expansion
 of \cite{Bruno6} gives the asymptotics when $x\to 0$ far from the poles.

\vskip 0.3 cm 
{\bf 2)} In this paper, section
   \ref{FULL}, the  recursive procedure is given to compute at any order the 
   expansions, for $x\to$~critical point, of the  4 and 3 real parameters solutions (namely, solutions such that 
$
\hbox{Tr}(M_jM_k)\neq \pm 2$, $0\leq \Re \sigma_i \leq 1$, $\sigma_i\neq 0,1
$). 
 The ordering of the terms in the expansion is sensibly
 depending on the initial conditions (i.e. on the exponent of the
 leading term). For this reason, so  far 
it has been thought that the expansion is formally uncomputable in 
 general. It is shown that this is not the case. The procedure to
 compute it in general is given, independently on the initial conditions
 (i.e. the value of $\sigma_i$).  
The convergent expansions  for $x\to 0$ are:
$$
 \delta(x)=\sum_{n=0}^\infty x^n\sum_{m=-n}^{n+2}
 \tilde{c}_{nm}(\sigma,a,\alpha,\beta,\gamma,\delta)x^{m\sigma}-1=
O(\hbox{\rm max}\{x^{\Re\sigma},x^{1-\Re\sigma}\}),~~~\tilde{c}_{00}=1
$$
for the solution (\ref{INTROO1}).
$$
\delta^{*}(x)= 
\sum_{n=1}^\infty x^n \sum_{m=-n-1}^{n+1} b_{nm}(\sigma,\phi,\alpha,\beta,\gamma,\delta) x^{m\sigma}=O(x)
$$
for the solution (\ref{INTROO2}).
$$
\delta^{*}(x)=\sum_{n=1}^\infty x^n \sum_{m=-n-1}^{n+1}
d_{nm}(\sigma,\phi,\alpha,\beta,\gamma,\delta) ~
x^{m(1-\sigma)}=O(x)
$$
for the solution (\ref{INTROO3}). The procedure is given to compute the
coefficients 
$\tilde{c}_{nm},b_{nm},d_{nm}$ in section \ref{FULL}. They are rational
functions of the integration constants $\sigma$, $\exp\{i\phi\}$ and of $\alpha,\beta,\gamma,\delta$. 

\vskip 0.3 cm 
{\bf 3)}  In this paper are also computed   
 the explicit formulae which
express Tr$(M_jM_k)$ as functions of the coefficients of PVI and of
the integration constants for 4 and 3 real parameter solutions, namely for  $0\leq \Re\sigma_i\leq 1$, $ \sigma_i\neq 0,1$ (see
 (\ref{COPPIAA}) for $0\leq \Re\sigma_i < 1$ and (\ref{COPPIAA1}) for
 $\Re\sigma_i=1$).  This is the first time that the explicit formulae appear in
 the literature for the general PVI 
(for the special case   $\beta=\gamma=0$, $\delta=1/2$  they are
 given in 
\cite{DM} and \cite{D4}).   These formulae are necessary
 for the  solution of the connection problem.

\vskip 0.3 cm

  The relevance of this paper is that, together with  all previous contributions, first of all that of Jimbo \cite{Jimbo} and then  the series of  papers  \cite{D4} \cite{D3} \cite{D2} \cite{D1},  PVI may be considered "solved" (in the meaning stated in the introduction), solved in all the
cases when there is a one-to-one correspondence between monodromy data
and Painlev\'e transcendents and $<M_0,M_x,M_1>$ is irreducible. This is because  {\it all} the critical
behaviors have been obtained,  and   {\it almost all} the parametrizations of the integration contants in terms of monodromy data. 
 ``Almost'' means that  some  special
values of the monodromy data $\theta_\mu$ (to be introduced in
section \ref{REVIEW}, see (\ref{puzzapi})) are poles
of the connection formulae. It is possible to compute the formulae in these  special cases as well, with no conceptual
changes in the general scheme of \cite{D2} and \cite{Jimbo}. These
  very time-consuming  computations  will be done only when   
one gets specifically interested  
in some  special case. In  \cite{DM} \cite{D4}, 
 all the computations for the relevant special
 case of PVI associated to a Frobenius manifold are done.     

\vskip 0.3 cm 
 This paper is organized as follows: 

- Section \ref{REVIEW}:  review the isomonodromy deformation approach to PVI.  

- Section \ref{prima}:  statement of the critical behavior (\ref{INTROO3}), when
  $x\to 0$ [Proposition \ref{prop3}]. Review of  (\ref{INTROO1}) and (\ref{INTROO2}).

- Section \ref{frazione}: proof of   (\ref{INTROO3}), via a symmetry of PVI which
  transforms (\ref{INTROO2}) into (\ref{INTROO3}).

- Section \ref{conection}: The connection problem. All the 
 formulae relating monodromy data and integration
  constants are given for the small power type behaviors, the sine-type oscillatory behaviors and inverse sine-type oscillatory behaviors. 

- Section \ref{EAUTUNNO}:  example of the above connection formulae for
  PVI associated to a Frobenius manifold. 

- Section \ref{FULL}: recursive computation of the full expansion of $y(x)$ (of
  $\delta(x)$, $\delta^*(x)$) at the critical points. 

- Appendix I: review of the procedure of Jimbo to obtain 
(\ref{INTROO1}) and (\ref{INTROO2}). 

- Appendix II: review of the elliptic representation and proof of  the
  convergence of the full expansion of $y(x)$ (of
  $\delta(x)$, $\delta^*(x)$).

\section{Review of the Isomonodromy Deformations}   
\label{REVIEW}

(PVI) is the isomonodromy deformation equation of the $2\times 2$ matrix linear Fuchsian system of ODEs given in equation (\ref{SYSTEM}). The  $2\times 2$ matrices  $A_i(x)$  depend on the 
 parameters $\alpha,\beta,\gamma,\delta$ according to 
the following relations:  
$$
 A_0+A_1+A_x = -{\theta_{\infty}\over 2}
 \sigma_3,~~\theta_\infty\neq 0,
$$
$$
\hbox{ Eigenvalues}~( A_i) =\pm {1\over 2} \theta_i, ~~~i=0,1,x;
$$
\be
\label{puzzapi}
\theta_0^2=-2\beta,~~~\theta_x^2=1-2\delta,~~~\theta_1^2=2\gamma,~~~(\theta_\infty-1)^2=2\alpha
\ee
Here $\sigma_3:=$ diag$(1,-1)$ is the Pauli matrix. The condition $\theta_\infty\neq 0$ is not restrictive, because $\theta_\infty=0$ is equivalent to $\theta_\infty=2$.   
 The equations of monodromy preserving deformation (Schlesinger equations), can be written in Hamiltonian form and reduce
 to (PVI), being the transcendent $y(x)$ the solution $\lambda$  of
 $A(x,\lambda)_{1,2}=0$. Namely:
\be
y(x)= 
{x~(A_0)_{12} \over x~\left[
(A_0)_{12}+(A_1)_{12}
\right]- (A_1)_{12}},
\label{leadingtermaprile}
\ee
The matrices $A_i(x)$, $i=0,x,1$, 
 depend on $y(x)$, ${d y(x)\over dx}$ and $\int y(x)$ 
through rational functions, which are given in \cite{JMU}

\vskip 0.2 cm 
The standard choice of a fundamental matrix $\Psi$ is as 
follows: 
\be
\Psi(\lambda)= 
\left\{ 
\matrix{
\left[
I+O\left({1\over \lambda}\right)
\right]~\lambda^{-{\theta_\infty\over 2}\sigma_3} \lambda^{R_\infty},&~~~\lambda\to\infty;
\cr
\cr
\psi_0(x) \bigl[I+O(\lambda)\bigr]~\lambda^{{\theta_0\over
    2}\sigma_3}\lambda^{R_0}C_0,&~~~\lambda\to 0;
\cr
\cr
\psi_x(x)\bigl[I+O(\lambda-x)\bigr]~(\lambda-x)^{{\theta_x\over
    2}\sigma_3}(\lambda-x)^{R_x}C_x,&~~~\lambda\to x;
\cr
\cr
\psi_1(x)\bigl[I+O(\lambda-1)\bigr]~(\lambda-1)^{{\theta_1\over
    2}\sigma_3}(\lambda-1)^{R_1}C_1,&~~~\lambda\to 1;
}\right.
\label{PSIlocalelo}
\ee
Here $\psi_0(x)$,  $\psi_x(x)$, $\psi_1(x)$ are the $2\times 2$ diagonalizing
matrices of $A_0(x)$, $A_1(x)$, $A_x(x)$ respectively. They are
defined by multiplication to the right by arbitrary diagonal matrices,
possibly depending on $x$. $C_\nu$, $\nu=\infty,0,x,1$,
 are invertible {\it connection matrices}, independent of $x$
\cite{JMU}. 
Each $R_\nu$,  
$\nu=\infty,0,x,1$, is also independent of $x$, and:
$$
R_\nu=0 \hbox{ if } \theta_\nu\not\in {\bf Z},~~~~~
R_\nu=\left\{
\matrix{
\pmatrix{0 & *\cr 0 & 0},~~~ \hbox{ if } \theta_\nu>0 \hbox{
  integer}
 \cr
\cr
\pmatrix{0 & 0\cr * & 0},~~~ \hbox{ if } \theta_\nu<0 \hbox{
  integer}
}
\right.
$$
If $\theta_i=0$, $i=0,x,1$, then 
 $R_i$ is to be considered the Jordan form $\pmatrix{0 & 1 \cr 0 & 0}$
of $A_i$.

\begin{figure}
\epsfxsize=12cm

\epsfysize=12cm

\centerline{\epsffile{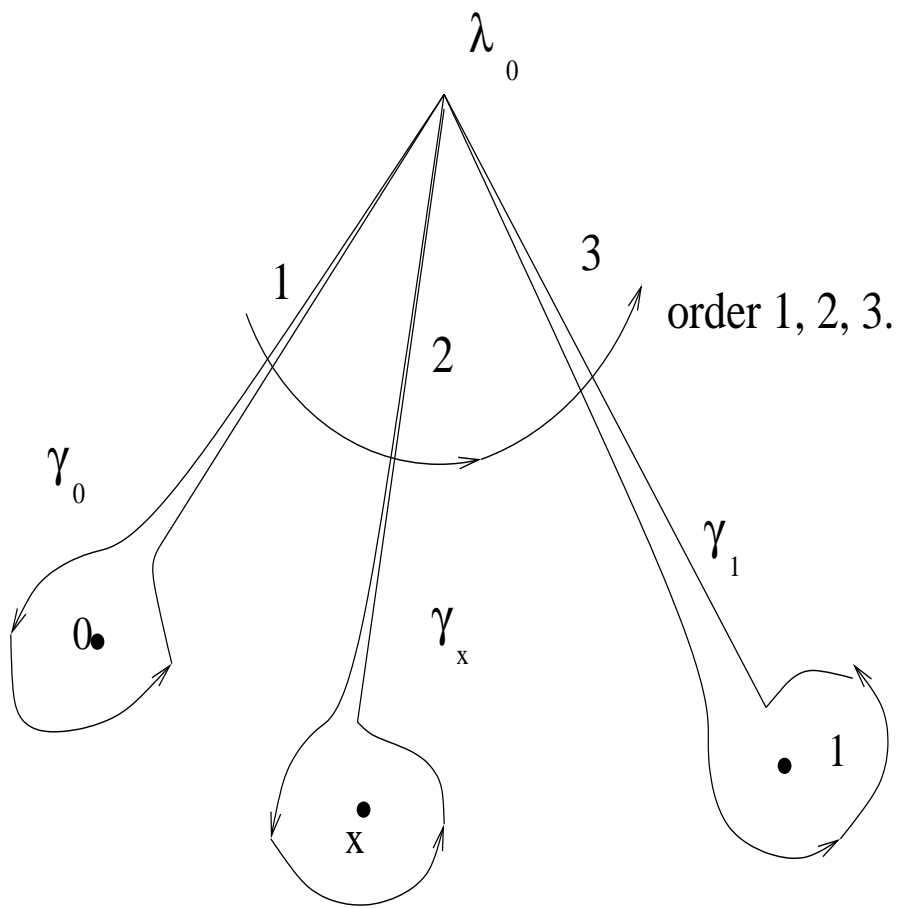}}
\caption{The ordered basis of loops $\Gamma$}
\label{figure1}
\end{figure}

\begin{figure}
\epsfxsize=12cm

\epsfysize=12cm

\centerline{\epsffile{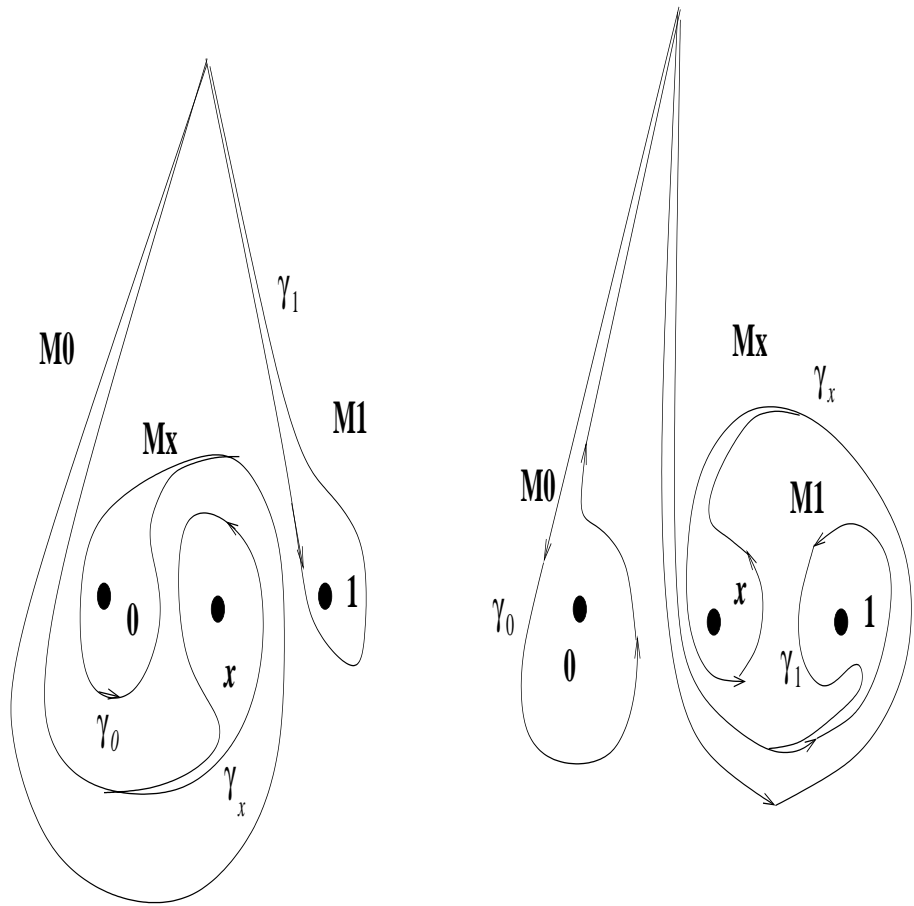}}
\caption{Other choices of the ordered basis of loops, $\Gamma_0$ (left) and $\Gamma_1$ (right)}
\label{figure2}
\end{figure}

\vskip 0.2 cm 
 Let a basis of loops in the order $(1,2,3)=(0,x,1)$ be fixed. There are several (infinite)  choices of such a basis. Here the basis $\Gamma$  of figure
 \ref{figure1} is chosen (other possible simple choices are the basis $\Gamma_0$ and $\Gamma_1$in figure \ref{figure2}). 

  Let the $x$-plane be cut by the condition that $|\arg x|<\pi$, $|\arg (1-x)|<\pi$, so that $A(x,\lambda)$ and $y(x)$ make sense as $x$-branches.
 
 When $\lambda$ goes around a counter-clockwise loop around $0,x,1$,
 then $\Psi$ is transformed by right multiplication by {\it monodromy
 matrices} $M_0, M_x, M_1$:
$$ 
\Psi\mapsto \Psi M_j,~~~M_j=C_j^{-1} \exp\{i\pi\theta_j\sigma_3\}\exp\{2\pi i R_j\} C_j,~~~j=0,x,1.
$$
For the loop $\gamma_\infty$: $\lambda \mapsto \lambda e^{-2\pi i}$,
$|\lambda|>\max\{1,|x|\}$,  the monodromy at infinity is:  
$$
M_\infty=\exp\{i\pi\theta_\infty\}~\exp\{- 2\pi i R_\infty\}. 
$$ 
The following relation holds:
$$
\gamma_0\gamma_x\gamma_1\gamma_\infty=1,~~~M_1M_x M_0 M_\infty=I
$$


 \vskip 0.2 cm 
\noindent
 The {\bf monodromy data} of the fuchsian system, {\it with respect to a basis of loops} $\Gamma$,  are the following set
of data:
\vskip 0.2 cm

a) The exponents  $\pm\theta_0,\pm\theta_x,\pm\theta_1,\pm(\theta_{\infty}-1)$, with  
$\theta_\infty \neq 0$. 

b) Matrices  $R_0,R_x,R_1,
R_\infty$, such that:  
$$
R_\nu=0 \hbox{ if } \theta_\nu\not\in {\bf Z},~~~~~
R_\nu=\left\{
\matrix{
\pmatrix{0 & *\cr 0 & 0},~~~ \hbox{ if } \theta_\nu>0 \hbox{
  integer}
 \cr
\cr
\pmatrix{0 & 0\cr * & 0},~~~ \hbox{ if } \theta_\nu<0 \hbox{
  integer}
}
\right.
$$
$$
R_j=\pmatrix{0 & 1 \cr 0 & 0},~~~\hbox{ if } \theta_j=0,~~~j=0,x,1.
$$

c) three  monodromy matrices $M_0$, $M_x$, $M_1$ relative to the 
loops, similar to the matrices $\exp\{ i\pi\theta_i\sigma_3\}
\exp\{ 2\pi i R_i\}$, $i=0,x,1$, satisfying:  
$$ 
  M_1~M_x~M_0 = e^{-i\pi \theta_{\infty}\sigma_3 }e^{2 \pi i  R_{\infty}}
$$

The data 
$\pm\theta_0,\pm\theta_x,\pm\theta_1,\pm(\theta_\infty-1)$ are fixed by
the equation. The other monodromy data are free. To each choice of
them, there corresponds at least one  fuchsian system (the solution of a
Riemann-Hilbert problem for the given monodromy data). To such a 
fuchsian system,  a branch $y(x)$ is associated. 
Therefore, there is a correspondence between 
a set of monodromy data and a branch 
$y(x)$. In generic cases, the correspondence is 
one-to-one. This is stated in the following theorem,
 proved in \cite{D1}, section 3: 

Let $\Theta,R,M$ stands for the collection  $\theta_0,\theta_x,\theta_1,\theta_\infty\neq 0$, $R_0,R_x,R_1,R_\infty$, $M_0,M_x,M_1$. 

\bth
{\it Let a basis of loops $\Gamma$  be chosen and let the monodromy data with respect to $\Gamma$ be   $\Theta,R,M$ 
  satisfying a), b), c) above. There is a one to one correspondence between the monodromy data 
 and one branch of a transcendent $y(x)$,  except  when 
at least one $\theta_\nu\in{\bf Z}\backslash\{0\}$ and simultaneously  
$R_\nu=0$.}
The branch in one to one correspondence with $\Gamma$,  $\Theta,R,M$  will be denoted:
\be
\label{abranch0}
 y(x)=f_{\Gamma}(x;\Theta,R,M)
\ee
\label{pro1}
\eth

\vskip 0.2 cm

Note that for $\theta_j=0$, $M_j$ can be put in Jordan form $\pmatrix{1 & 2\pi i \cr
0 & 1}$. Therefore: 

\noindent
{\it There is a  one to one correspondence 
if and only if one of the following conditions is satisfied: 
\vskip 0.2 cm 
(1) $\theta_\nu \not\in {\bf Z}$,  for every  $\nu=0,x,1,\infty$;
\vskip 0.2 cm 
(2)  some $\theta_\nu\in {\bf Z}$ and  $R_\nu\neq 0$, $\theta_\nu\neq 0$
\vskip 0.2 cm 
(3)  some $\theta_j=0$ ($j=0,x,1$) and simultaneously $\theta_\infty\not \in {\bf Z}$, or $\theta_\infty\in{\bf Z}$ and $R_\infty\neq 0$. 
}

\vskip 0.2 cm 
\noindent 
Equivalently: {\it There is one to one correspondence except when one of the 
matrices $M_i$ ($i=0,x,1$), or $M_\infty^{-1}=M_1M_xM_0$, is equal to $\pm I$. 
}

\vskip 0.5 cm
\noindent 
Define the following quantities:  
\be
\label{UFFA!}
p_\mu= \hbox{Tr} M_\mu=2\cos(\pi\theta_\mu),~~~p_{ij}=\hbox{Tr}(M_iM_j),~~~ 
 \mu=0,x,1,\infty,~~~ i,j\in \{0,x,1\}
\ee
These  coordinates  describe the space of monodromy  data,
which is an affine cubic surface \cite{Iwa} \cite{Jimbo}:
$$
p_{0x}^2+p_{01}^2+p_{x1}^2+p_{0x}p_{01}p_{x1}-(p_0p_x+p_1p_\infty)p_{0x}-(p_0p_1+
p_xp_\infty)p_{01}-(p_xp_1+p_0p_\infty)p_{x1}+
$$
$$
+p_0^2+p_1^2+p_x^2+p_\infty^2+p_0p_xp_1p_\infty -4 =0
$$
The above relation follows by taking the trace of the relation $M_1M_xM_0M_\infty=I$.  
\vskip 0.2 cm 
{\it If the monodromy group $<M_0,M_x,M_1>$ is not reducible, or one of the matrices $M_0,M_x,M_1,M_1M_x M_0$ is not the identity, the above $p_\mu$'s, $p_{ij}$'s are a good system of coordinates for the monodromy group} \cite{Iwa},\cite{Jimbo}.   

\vskip 0.2 cm 
As a consequence, a branch of a transcendent is uniquely parametrized by
the $p_\mu$'s and $p_{ij}$'s to which it is in one to one correspondence. In other words, the integration constants are
 functions of $p_\mu$'s (or
$\theta_\mu$'s)  and $p_{ij}$'s. The following notation expresses this dependence\footnote{The integration constants are two complex numbers. $ \theta_0,\theta_x,\theta_1,\theta_\infty$ are fixed by the equation and $p_{0x},p_{01},p_{x1}$ are not independent, becaus of the cubic surface relation. Accordingly, only two complex parameters are free.}:
\be
\label{abranch}
 y(x)=f_\Gamma(x;\theta_0,\theta_x,\theta_1,\theta_\infty,p_{0x},p_{01},p_{x1})
\ee
 A remarkable fact, established in  Jimbo's paper \cite{Jimbo}, is that {\it this parametrization is explicit}, namely the integration constants are elementary or classical  transcendental functions of the monodromy data. 

As a consequence of this explicit parametrization of the three couples of
 integration constants at the three critical points in terms  of {\it the same} monodromy data, the connection problem is solved.   This is precisely the power of the method of
 monodromy preserving deformations. 

\vskip 0.2 cm 
We observe that, when the monodromy group is reducible, but none of the monodromy matrices $M_0,M_x,M_1,M_1M_x M_0$ is the identity, the one to one correspondence still holds, but the $p_\mu$'s, $p_{ij}$'s are no longer a good parametrization. The solutions in this case are known (see the Riccati solutions \cite{watanabe}, \cite{mazzoccoratio}). 

\subsubsection{Analytic continuation of a branch}

 It is to be stressed   that (\ref{abranch0}) (or (\ref{abranch})) represents {\it a branch} of a transcendent, for $|\arg x|<\pi$, $|\arg (1-x)|<\pi$, uniquely identified by the parametrization in terms of monodromy data $\Theta,R,M$, which are associated to a the basis $\Gamma$.  We show below that if the same monodromy data $\Theta,R,M$ are associated to another basis   basis $\Gamma^\prime$, like $\Gamma_0$ and $\Gamma_1$ in figure \ref{figure2},  the parametrization $f_{\Gamma^\prime}(x;\Theta,R,M)$ is the branch  $y^{\prime}(x)$  of the analytic continuation of $y(x)$  along a path in the $x$ plane. Such path is the path that induces  in the $\lambda$-plane (as $x$ moves 
in the $\lambda$-plane around $\lambda=0$ or 1)  the deformation of  the basis $\Gamma$ into $\Gamma^\prime$. 

 The two basis $\Gamma_0$ on the left and $\Gamma_1$ on the right of  figure \ref{figure2} can be regarded as the deformation of $\Gamma$,  when $x$ goes around a counterclockwise loop around $\lambda=0$ or 1 respectively, namely when  $x$ goes counterclockwise along a path around $x=0$ or $x=1$ in the $x$-plane. 
The branch (\ref{abranch0})  undergoes its analytic continuation along these paths. {\it Being the deformation isomonodromic}, the monodromy matrices after the deformation do not change. So, the same $M_0,M_x,M_1$ are also  assigned to $\gamma_0,\gamma_x,\gamma_1$ of the basis $\Gamma_0$ or $\Gamma_1$.  
Let $\tilde{y}(\tilde{x})$ represent the analytic continuation of $y(x)$ of (\ref{abranch0}). It is defined on the universal covering of points $\tilde{x}$ and can be written as $\tilde{y}(\tilde{x})=f_\Gamma(\tilde{x};\Theta,R,M)$. Its branch $y^\prime(x)$  for $|\arg x|<\pi$, $|\arg (1-x)|<\pi$,  has again a parametrization in terms of of $M_0,M_x,M_1$. But it differs form  (\ref{abranch0}), because it is computed w.r.t. the basis $\Gamma_0$ or $\Gamma_1$. Let it be denoted by:  
$$ 
y^\prime(x)=f_{\Gamma^\prime}(x,\Theta,R,M)~~~~\hbox{ where $\Gamma^\prime$ stands for $\Gamma_0$ or $\Gamma_1$}.
$$
One has now to compute $f_{\Gamma^\prime}(x,\Theta,R,M)$.  The way to do this is to compute the monodromy  matrices associated to the  basis  $\Gamma$, being   $M_0,M_x,M_1$  associated to $ \Gamma_0$ or $\Gamma_1$.   In order to do this, observe that  the loops of figure \ref{figure1} can be written as a product of the loops of figure \ref{figure2} as:
$$ 
\gamma_0(\hbox{ of figure \ref{figure1}})= \gamma_x^{-1}\gamma_0\gamma_x( \hbox{ of figure \ref{figure2}, left}),
$$
$$ 
\gamma_x(\hbox{ of figure \ref{figure1}})= \gamma_x^{-1}\gamma_0^{-1}\gamma_x \gamma_0\gamma_x( \hbox{ of  figure \ref{figure2}, left}),
$$
$$ 
\gamma_1(\hbox{ of figure \ref{figure1}})= \gamma_1( \hbox{ of left figure \ref{figure2}, left}),
$$
and:
$$ 
\gamma_0(\hbox{ of figure \ref{figure1}})= \gamma_0( \hbox{ of figure \ref{figure2}, right}),
$$
$$ 
\gamma_x(\hbox{ of figure \ref{figure1}})= \gamma_1^{-1}\gamma_x \gamma_1( \hbox{ of  figure \ref{figure2}, right}),
$$
$$ 
\gamma_1(\hbox{ of figure \ref{figure1}})= \gamma_1^{-1}\gamma_x^{-1}\gamma_1 \gamma_x \gamma_1( \hbox{ of left figure \ref{figure2}, right})
$$
 It follows that, beeing $M_0,M_xM_1$ the monodromy matrices for the basis $\Gamma_0$ or $\Gamma_1$,  the monodromy matrices with respect to the initial basis $\Gamma$ are:
\be
\label{ac1}
M_0^\prime=M_xM_0M_x^{-1},~~~M_x^\prime=M_xM_0M_xM_0^{-1}M_x^{-1},~~~M_1^\prime=M_1,
\ee
in the left case (basis $\Gamma_0$, counterclockwise loop of $x$ around 0). 
 \be
\label{ac2}
M_0^\prime=M_0,~~~M_x^\prime=M_1M_x M_1^{-1},~~~M_1^\prime = M_1 M_x M_1 M_x^{-1} M_1^{-1}, 
 \ee
in the right case (basis $\Gamma_1$, counterclockwise loop of $x$ around 1). The above transformation of the monodromy matrices is an action of the {\it braid group}. It implies that $\Theta$ and $R$ are not changed. The branch of the analytic continuation is then:
$$
y^\prime(x)=f_{\Gamma}(x;\Theta,R,M^\prime)
$$
and thus the computation of  $f_{\Gamma^\prime}(x;\Theta,R,M)$ has been completed. To summarize:

\vskip 0.2 cm 
{\it Let $\Theta,R,M$ be given. The choice of the basis $\Gamma$ determines a branch $y(x)=f_{\Gamma}(x;\Theta,R,M)$.  The choice of another basis $\Gamma^\prime$ determines another branch $y^\prime(x)=f_{\Gamma^\prime}(x;\Theta,R,M)$, which is a branch of the analytic continuation of $y(x)$ along the path of $x$ which deforms $\Gamma$ to $\Gamma^\prime$. The relation is 
$$
  f_{\Gamma^\prime}(x;\Theta,R,M) =f_{\Gamma}(x;\Theta,R,M^\prime) 
$$
where $M\mapsto M^\prime$ is an action of the braid group generated by (\ref{ac1}) and (\ref{ac2}).
In other words, the anayltic continuation of $ y(x)=f_{\Gamma}(x;\Theta,R,M)$ is $y^\prime(x)=f_{\Gamma}(x;\Theta,R,M^\prime)$.
}

\vskip 0.2 cm 
In terms of the coordinates $p_\mu$'s and $p_{ij}$"s the above transformation of the matrices reads:  
 \be
\label{vesc}
\left\{
\matrix{
p_{x1}^\prime =
p_{x1}(p_{0x}^2-1)+p_{0x}p_{01}-(p_\infty p_x +p_1p_0)p_{0x}+p_\infty
p_0+ p_1 p_x
\cr
\cr
 p_{0x}^\prime =p_{0x},~~~
p_{01}^\prime = -p_{01}-p_{x1}p_{0x} +p_\infty p_x+p_1 p_0
}
\right.
\ee in the left case  (basis $\Gamma_0$,  $x$ goes around
 a loop around $0$), and:   
\be
\label{vesc1}
\left\{
\matrix{
p_{01}^\prime = 
p_{01}(p_{x1}^2-1)+p_{0x}p_{x1}-(p_\infty p_1 +p_0p_x)p_{x1}+p_\infty
p_x +p_0 p_1
\cr
\cr
 p_{1x}^\prime = p_{1x},
~~~p_{0x}^\prime =  -p_{0x}-p_{01}p_{x1} +p_\infty
 p_1 +p_0 p_x
}
\right.
\ee
in the right case (basis $\Gamma_1$, $x$ goes around
 a loop around $1$). The branch of the analytic continuation has parametrization: 
$$
y^\prime(x)=f_{\Gamma}(x;\theta_0,\theta_x,\theta_1,\theta_\infty,p_{0x}^\prime,p_{01}^\prime,p_{x1}^\prime),~~~~~|\arg x|<\pi,~~~|\arg(1-x)|<\pi.
$$
 The parametrization of the branch of the analytic continuation along more complicated paths is given by a suitable composition of (\ref{vesc}) and (\ref{vesc1}).

\vskip 0.2 cm 
As a final remark it is to be noted that the choice $|\arg x+2\pi k|<\pi$, $|\arg(1-x)+2\pi l|<\pi$, for some $k,l\in{\bf Z}$, is also possible. But the computation of the {\it explicit}  parametrization is done by the procedure of \cite{Jimbo}, which makes use of a reduction of (\ref{SYSTEM}) to hyper-geometric equations  (and by its  generalization of \cite{D2} to non-hyper-geometric reductions in case of Taylor solutions). This computation requires  that $k=l=0$. Accordingly, the formulas which parametrize the critical behaviors  in this paper are given for the branches $|\arg x|<\pi$, $|\arg(1-x)|<\pi$.

\section{Critical behavior at $x=0$}\label{prima}

In the following, it is understood that   $x \to$ critical point 
inside a sector. 
The behavior of $y(x)$ at   $x=0,1,\infty$
 is determined by three {\it critical exponents}
 $\sigma_0,\sigma_1,\sigma_\infty$ respectively,  given by: 
$$
 2\cos(\pi \sigma_0)= p_{0x},
~~~ 2\cos(\pi \sigma_1)= p_{x1},~~~
 2\cos(\pi \sigma_\infty)=p_{01},~~~0\leq\Re\sigma_i\leq 1, 
$$
where $p_{ij}$ are (\ref{UFFA!}). 

\vskip 0.2 cm 
\noindent
{\bf Remark: }
 The above relation
determines $\sigma_i$ up to $\sigma_i\mapsto \pm \sigma_i +2n$, $n\in
{\bf Z}$. One can then restrict to the case $0\leq \Re\sigma_i \leq 1$,
 as it is explained in \cite{D4} \cite{D3}. 
Despite this condition,  when $\Re\sigma_i=0$ the ambiguity of sign  cannot be eliminated. Namely:  
$$
 \sigma_i= \pm i\nu,~~~\nu\in{\bf R}, ~~~~~~ p_{kl}=
 \cosh \pi \nu > 2,
$$
In case $\Re \sigma_i=1$ the ambiguity $\sigma_i \mapsto 2-\sigma_i$ cannot
be eliminated. Namely: 
$$ 
 \sigma_i= 1\pm i\nu,~~~\nu\in{\bf R}, ~~~~~~ p_{kl}=
 -\cosh \pi \nu <-2,
$$
Anyway, a solution $y(x)$ corresponding to such monodromy data is
invariant for the change of sign of $\nu$, as it will be explained below.

\vskip 0.3 cm 
We start with the critical point $x=0$. In the following, we use the notation
$\sigma:=\sigma_0$.   Let also $|x|<\epsilon<1$, where $\epsilon$  is 
   sufficiently small for all our purposes.  The first result of this paper is the following:

\bpr
\label{prop3} {\rm[{\bf Inverse sine-type oscillatory behaviors}]} The equation PVI admits solutions with a branch at $x=0$  behaving in the
following way  when 
 $x\to 0$, with arg$(x)$ bounded:  
\be
\label{LABEL}
y(x)= {1\over iA \sin\Bigl(i(1-\sigma)\ln x+\phi\Bigr)+B
  +\delta^{*}(x)},~~~\delta^*(x) =O(x)
\ee
where $\sigma,\phi\in{\bf C}$ are the integration constants, satisfying
  $\Re \sigma=1$, $\sigma\neq 1$. The coefficients $A$ and $B$ are: 
$$
B 
  =
     {\nu^2+2\gamma-2\alpha\over 2
	\nu^2}
=  {\nu^2+(\theta_1)^2-(\theta_{\infty}-1)^2\over 2
	\nu^2}
,
~~~~~\sigma = 1+ i \nu ,~~~\nu\in{\bf R},~~~\nu\neq 0,
$$

\be
\label{AA}
A = i\sqrt{{2\alpha\over \nu^2}+B^2} = 
i\sqrt{
{(\theta_\infty -1)^2\over \nu^2}
+
B^2
}=
\ee
$$
 = 
{\sqrt{
\Bigl[(1-\sigma)^2-(\theta_\infty-1-\theta_1)^2\Bigr]\Bigl[(\theta_\infty-1+
\theta_1)^2-(1-\sigma)^2\Bigr]}
\over
 2(1-\sigma)^2}
$$
The vanishing term $\delta^*(x)$ has convergent  expansion for
$0<|x|<\epsilon$: 
\be
\label{ddelta3}
{
\delta^{*}(x)=\sum_{n=1}^\infty x^n \sum_{m=-n-1}^{n+1} d_{nm}~
x^{m(1-\sigma)}
~ =
\sum_{m_1=1}^\infty \sum_{m_2=-1}^{2m_1+1}
   e_{m_1m_2}x^{m_1\sigma}x^{m_2(1-\sigma)}
}
\ee
$$
e_{m_1m_2}=d_{m_1,m_2-m_1}. 
$$
The coefficients  are certain rational
 functions of $\sigma$ and $\exp\{i\phi\}$,  which can be computed by
 direct substitution into PVI (see section \ref{FULL}). 
The constant $\sigma$ is related to the monodromy data associated to
$y(x)$ by: $2\cos(\pi\sigma)=p_{0x}<-2$.
\epr

\noindent
Since $\sin (2x)= 1-2 \sin^2(x-\pi/4)$, we can also rewrite:
\be
\label{XXX}
{y(x)=
\left\{
-2iA \sin^2\left(i{1-\sigma \over 2}\ln x +{\phi\over 2}-{\pi \over
  4}\right)~+iA+B~+\delta^{*}(x)
\right\}^{-1}
}
\ee
Let $r\in {\bf C}$, $r\neq 0$. It is convenient, for future developments,  to re-parametrize  $\phi$ in terms of
$r$ as follows (at this stage of the discussion, this may be
temporarily taken as the definition of $r$):
$$
\phi=i\ln{2r\over (1-\sigma)A}
$$
The reason to introduce $r$ is that it is a natural parameter that
  will be written  in section 
\ref{conection} as a function of the monodromy data  associated to the basis $\Gamma$  of figure \ref{figure1}: 
$$
 r=r(\sigma,\theta_0,\theta_x,\theta_1,\theta_\infty,p_{x1},p_{01}).
$$ 
This parametrization  identifies uniquely the branch.  

The sign of the square root $A$ can be chosen arbitrarily, because it
changes $\phi \mapsto \phi +(2k+1)\pi$, and $y(x)$ is invariant. 
It is to be noted that the condition $\Re \sigma=1$ does not fix  the ambiguity  $\sigma=1+i\nu\mapsto 1-i\nu$,
(namely $\sigma\mapsto 2-\sigma$), but the substitution $\sigma\mapsto 2-\sigma$ induces $\phi \mapsto -\phi +(2k+1)\pi$ (see Appendix I, subsection \ref{maldipiedi}), and thus $y(x)$ is invariant.

\vskip 0.3 cm

\noindent
{\bf Remark:} {\it We must  keep $\delta(x)$ in the denominator}. 
This term is essential in
that it determines the position of the movable poles, which occur when
 the
denominator vanishes at some isolated points.

\vskip 0.3 cm 
For completeness, the  results about the
critical behavior when $0\leq
\Re\sigma<1$, $\sigma\neq 0$, are reported below. Though the critical behaviors are already known and appear in 
 \cite{Jimbo}, \cite{DM} \cite{D4} \cite{D3}
\cite{D2},  
the  expansions of the terms $\delta(x)$ and $\delta^*(x)$
 in the propositions 
below is a result of the present paper (see section \ref{FULL}).

\bpr
\label{prop1} {\rm [{\bf Small power type behaviors (Jimbo)}]}  The equation PVI admits solutions with a branch having the following
behavior, when $x\to 0$, arg$(x)$ bounded (\cite{Jimbo}, 
\cite{DM} \cite{D4} \cite{D3} \cite{D2}): 
\be
\label{giusta}
{y(x)= a x^{1-\sigma} ~(
1+ \delta(x)
)},~~~\delta(x) =O(\hbox{\rm max}\{x^{1-\Re
\sigma},x^{\Re\sigma}\})
\ee
where $ a, \sigma \in {\bf C}$ are integration constants such that
$a\neq 0$ and $0<\Re \sigma<1$. The higher order term 
$\delta(x)$ has the following 
convergent expansion for $0<|x|<\epsilon$ (section \ref{FULL}): 
$$
{
\delta(x)=-1+\sum_{n=0}^\infty x^n\sum_{m=-n}^{n+2}
 \tilde{c}_{nm}x^{m\sigma},~~~\tilde{c}_{00}=1}
$$
We can also write: 
\be
\label{ddelta1}
{
\delta(x) =  \sum_{m_2=0}^\infty \sum_{m_1=0}^{2m_2+2}
\delta_{m_1m_2}x^{m_1\sigma}x^{m_2(1-\sigma)}
 ,~~~~m_1+m_2\geq 1}
\ee
$$
\delta_{m_1m_2}=\tilde{c}_{m_2,m_1-m_2}.
$$ 
The coefficients  are certain rational
 functions of $\sigma$ and $a$,  which can be computed by
 direct substitution into PVI (see section \ref{FULL}). 
The exponent $\sigma$ is related to the monodromy data associated to
$y(x)$ by: $2\cos(\pi\sigma)=p_{0x}$.
\epr

As before, we  re-parameterize $a$ in terms of a new $r\in {\bf C}$:
$$
a={
1
\over
16\sigma^3r
}\Bigl[\sigma^2-(\sqrt{-2\beta}-\sqrt{1-2\delta}~)^2\Bigr]\Bigl[(\sqrt{-2\beta}+
\sqrt{1-2\delta}~)^2-\sigma^2\Bigr]
=
$$
\be
\label{X}
={[\sigma^2-(\theta_0-\theta_x)^2][(\theta_0+\theta_x)^2-\sigma^2]
\over
16\sigma^3 r}.
\ee
$r$ will be naturally introduced when proving (\ref{giusta}) in
Appendix I.  The parametrization of $r$ in terms of monodromy data identifies the branch uniquely.

\vskip 0.2 cm 
\noindent
{\bf Remark:} For special values of $\sigma$ we have the following solutions:
\be
y(x)={\theta_0\over \theta_0+\theta_x}~x ~\mp~{r\over\theta_0+\theta_x}
~x^{1+\sigma}+O(x^2),~~~\sigma=\pm(\theta_0+\theta_x)\neq 0,
\label{passero}
\ee
\be
y(x) = {\theta_0\over \theta_0-\theta_x}~x~\mp
~{r\over\theta_0-\theta_x}~x^{1+\sigma}+O(x^2),
~~~
\sigma=\pm(\theta_0-\theta_x)\neq 0.
\label{aquila}
\ee

\vskip 0.3 cm

\bpr
\label{prop2} {\rm [{\bf Sine-type oscillatory behaviors}]} 
  The equation PVI admits solutions with a branch having  the following
behavior, when $x\to 0$, arg$(x)$ bounded (\cite{Jimbo},  \cite{D2}): 
\be
\label{mifido}
 y(x) = x ~\Bigl\{ 
iA\sin(i\sigma \ln x + \phi) + B + \delta^*(x)
\Bigr\},~~~\delta^*(x)=O(x)
\ee
where $\sigma,\phi\in{\bf C}$ are integration constants such that $\Re
\sigma =0$, $\sigma\neq 0$. The coefficients are: 
$$
B = 
{\theta_0^2-\theta_x^2 + \sigma^2 \over 2\sigma^2}    
 ={\sigma^2-2\beta-1+2\delta \over 2\sigma^2} 
,
$$

$$
A={ \sqrt{
\bigr[\sigma^2-(\theta_0+\theta_x)^2
\bigl]
\bigr[
(\theta_0-\theta_x)^2-\sigma^2
\bigl]
}
\over 2\sigma^2
}=
 \sqrt{ {\theta_0^2\over \sigma^2} -B^2 }=\sqrt{- {2\beta\over \sigma^2} -B^2 },
$$
The term $\delta^*(x)$ has convergent expansion for $0<|x|<\epsilon$: 
\be
\label{ddelta2}
{
\delta^{*}(x)= 
\sum_{n=1}^\infty x^n \sum_{m=-n-1}^{n+1} b_{nm} x^{m\sigma}
=
\sum_{m_2=1}^\infty \sum_{m_1=-1}^{2m_2+1} a_{m_1m_2}
x^{m_1\sigma}x^{m_2(1-\sigma)},
}
\ee
$$
a_{m_1m_2}=b_{m_2,m_1-m_2}.
$$
The constant $\sigma$ is related to the monodromy data associated to
$y(x)$ by: $2\cos(\pi\sigma)=p_{0x}>2$.
\epr

\vskip 0.3 cm 
The critical  behavior can be also written as:
\be
\label{mmfido}
{y(x)= x\left\{
-2iA\sin^{2}\left(
i{\sigma\over 2}\ln x +{\phi\over 2}-{\pi \over 4}
\right)~+iA+B+\delta^{*}(x)
\right\}}
\ee

\vskip 0.2 cm 
 We rewrite $\phi$ in terms of the new integration constant $r$, which
 will be expressed as
 $r=r(\sigma,\theta_0,\theta_x,\theta_1,\theta_\infty,p_{x1},p_{01})$
 in section \ref{conection} : 
$$
\phi= i\ln {2r\over \sigma A}
$$
Any sign of the square root in $A$ can be chosen  (change $\phi \mapsto
\phi +(2k+1)\pi$). The condition $\Re \sigma =0$ does not fix the ambiguity  $\sigma \mapsto -\sigma$. Nevertheless, (\ref{mifido}) is  invariant for
$\sigma\mapsto -\sigma$, because this   induces the change  $\phi\mapsto -\phi
+(2k+1)\pi$, $k\in{\bf Z}$. See the  Appendix I (at the end of subsection \ref{maldipiedi}) for the proof.

\vskip 0.2 cm
\noindent
{\bf Remark:} 
If $\Re\sigma=0$ ($\sigma=i\nu$, $\nu\in{\bf R}$) but  $\sigma\in\{\theta_0+\theta_x, \theta_0-\theta_x, -\theta_0+\theta_x, -\theta_0-\theta_x\}$, the solution of  proposition \ref{prop2} becomes: 
\be
y(x)=x\left({\theta_0\over i\nu}-{r\over i\nu} x^{i\nu}\right) +O(x^2),~~~~~\sigma=\theta_0\pm\theta_x,
\label{passero1}
\ee
\be
 y(x)=x\left(-{\theta_0\over i\nu}-{r\over i\nu} x^{i\nu}\right) +O(x^2),~~~~~\sigma=-(\theta_0\pm\theta_x)
\label{aquila1}
\ee

\vskip 0.5 cm

\noindent
{\bf How we Prove the above propositions:} 

(\ref{giusta}) and (\ref{mifido}) are proved (though not explicitly
written) in \cite{Jimbo}. The proof is reviewed  in Appendix I.

Formula (\ref{LABEL}) is proved in section \ref{frazione}, where we show that it is the
image of (\ref{mmfido}) via a fractional linear transformation (\ref{infi}).

 In section \ref{FULL}   the recursive procedure is given to compute the full expansion of
 $y(x)$, and thus  the series (\ref{ddelta3}), (\ref{ddelta1}), (\ref{ddelta2}). Their convergence follows from the elliptic
 representation.  The elliptic representation of PVI is analytically studied in \cite{D3}. All the critical behaviors  of $y(x)$ are computed for $0\leq \Re \sigma
 \leq 1$, $\sigma\neq 0,1$. The convergence of the full expansions is proved. 

-- When $0<\Re \sigma <1$, the critical behavior and the full convergent
 expansion obtained from the
 elliptic representation ((\ref{AAAa}) in Appendix II),  
coincides with (\ref{giusta}). This  proves
 the convergence of  (\ref{ddelta1}). 

-- When $\Re \sigma =0$ and $\Re \sigma =1$, the critical behaviors 
computed in \cite{D3} depend on two integration constants $\sigma$,
$\phi_E$ (three real constants). They are (see Appendix II): 
$$
y(x)= x\left[\sin^2\left(i{\sigma\over 2}\ln x +\phi_E+
  \sum_{n\geq 1}c_n(\sigma)[e^{-2i\phi_E}x^\sigma]^n\right)+\delta^*_E(x)\right],
$$
$$
\Re\sigma=0,~~~|x|<\epsilon,~~~|e^{-2i\phi_E}x^\sigma|<\epsilon,
~~~\delta^*_E(x)=O(x).
$$

$$
y(x)= \left[\sin^2\left(i{1-\sigma\over 2}\ln x 
 +\phi_E+ 
\sum_{n\geq 1} c_n(\sigma) [e^{-2i\phi_E}x^{1-\sigma}]^n
\right)+\delta^*_E(x) \right]^{-1},
$$
$$
\Re\sigma=1,~~~|x|<\epsilon,~~~|e^{-2i\phi_E}x^{1-\sigma}|<\epsilon,
~~~\delta^*_E(x)=O(x).
$$
 The series in $\sin^2(~..~)$  are absolutely convergent for sufficiently
 small $r<1$.  They are oscillating series that do not vanish when
 $x\to0$. In Appendix II the convergent
 expansion of the terms $\delta^*_E(x)$ is also given.  
 In subsection  \ref{oscilla} of Appendix II, the reader finds the proofs that the above behaviors
 coincide with our (\ref{mmfido}) and (\ref{XXX}). 
In order to do this, first  write $\sigma= -i\nu$ or $1+i\nu$, $\nu\in{\bf
 R}$. Then, it is shown that:
$$
-2iA\sin^2\left(
{\nu\over 2} \ln x +{\phi\over 2}-{\pi\over 4}
\right)~+iA+B = \sin^2\left(
{\nu\over 2} \ln x +f(x)
\right)
$$
where $f(x)$ is an oscillating function: 
$$
f(x)=\sum_{n\geq 0} f_n x^{-i\nu x},~~~f_n\in{\bf C}.
$$
 The coincidence of the result of the present paper with that of the elliptic representation,
 together with the convergence of the expansions of $\delta^*_E(x)$, 
proves  the
convergence of  (\ref{ddelta3}) and (\ref{ddelta2}).


\section{Proof of Proposition \ref{prop3}. The critical behavior at $x=0$ when  $\Re \sigma =
 1$. A Fractional Linear Transformation}
\label{frazione}
We consider the following fractional linear transformation, studied in
\cite{DM1}: 
\be
\label{infi}
   \theta_0^\prime = \theta_\infty -1,~~~ 
   \theta_x^\prime = \theta_1,~~~
   \theta_1^\prime = \theta_x,~~~
   \theta_\infty^\prime = \theta_0+1;~~~~~y^\prime(x^\prime)= {x\over y(x)},~~~x^\prime=x.
\ee
$$ 
(p_0,p_x,p_1,p_{\infty};p_{0x},p_{01},p_{x1})\mapsto$$
$$
\mapsto
(p_0^\prime,p_x^\prime,p_1^\prime,p_{\infty}^\prime;
p_{0x}^\prime,p_{01}^\prime,
p_{x1}^\prime)=(-p_\infty,p_1,p_x,-p_0;-p_{0x},-p_{01},p_{x1})
.$$
This is a symmetry of PVI, namely $y(x)$ solves PVI with
coefficients $\theta_\mu$ if and only if $y^\prime$ solves PVI with
coefficients $\theta^\prime_\mu$. 
We are going to use this transformation to obtain the critical
behavior of a transcendent $y^\prime(x)$  with $\Re
\sigma^\prime=1$, 
$p^\prime_{0x}<-2$, from the behavior of a
transcendent $y(x)$ with $p_{0x}=-p^\prime_{0x}>2$, $\Re
\sigma=0$.

\vskip 0.3 cm 
\noindent
$\diamond$ 
First, we compute the relation between $\sigma^\prime$ and
$\sigma$. The relation  $p^\prime_{0x}= -p_{0x}$ implies:  
$$
  2\cos (\pi \sigma^\prime) = - 2 \cos(\pi \sigma) ~~\Longrightarrow~~
\sigma^\prime=\pm \sigma +(2k+1)\pi ,~~k\in{\bf Z}
$$
The conditions $0\leq \Re\sigma<1$, $0\leq \Re\sigma^\prime<1$ imply
that:
\be
\label{infiinfi}
 {\sigma^\prime = 1-\sigma}
\ee

\vskip 0.3 cm 
\noindent
 $\diamond$ We compute the solution $y^{\prime}(x)$ with $\Re
 \sigma^\prime=1$  from the solution $y(x)$
 with  $\Re \sigma=0$. We know the critical behavior of this solution
 from the Jimbo's procedure of Appendix I: 
$$
y(x)= x~\bigr\{
iA \sin(i\sigma \ln x +\phi) +B +\delta^*(x)
\bigl\},~~~\sigma = \pm i \nu,
$$

$$
\phi= i\ln{2 r\over \sigma
  A},~~~~~~~~~~B={\nu^2+\theta_x^2-\theta_0^2\over 2\nu^2} \equiv
      {\nu^2+(\theta_1^\prime)^2-(\theta_{\infty}^\prime-1)^2\over 2
	\nu^2}
$$

$$
 A= i\sqrt{ {\theta_0^2\over \nu^2}+\left[
{\nu^2+\theta_x^2-\theta_0^2\over 2\nu^2}
\right]^2} 
           \equiv 
i\sqrt{
{(\theta_\infty^\prime -1)^2\over \nu^2}
+
\left[
{\nu^2+(\theta_1^\prime)^2-(\theta_\infty^\prime-1)^2\over 2\nu^2}
\right]^2
}
$$

The solution $y^\prime= {x\over y}$ obtained by fractional linear
transf. from $y(x)$ is immediately computed: 
$$
{
y^\prime(x)= \Bigl\{iA \sin\Bigl(i(1-\sigma^\prime)\ln x+\phi\Bigr)+B 
+\delta^*(x)
\Bigr\}^{-1}},~~~
\sigma^\prime = 1\mp i \nu ,$$
$$
~~~\phi=i\ln{2r\over (1-\sigma^\prime)A},
~~~\delta^*(x)=O(x)
$$ 

In section \ref{FULL} we compute the full expansion of $y(x)$, which
proves that 
$\delta^*(x)$ has the form  (\ref{ddelta2}). As a result, the
expansion of  
$y^\prime(x)$ obtained  by
the fractional linear transformation  proves 
(\ref{ddelta3}) from (\ref{ddelta2}).


\subsection{ The  case  $\Re \sigma = 1$ associated
 to a
  Frobenius Manifold}
\label{appendicite}

PVI is associated to a Frobenius Manifold when
  $\theta_0=\theta_x=\theta_1=0$ \cite{Dub}. 
The result of the general case, when $\Re \sigma=1$, becomes: 
$$
B= {\nu^2-(\theta_\infty-1)^2\over 2\nu^2},~~~~A=\pm i
{\nu^2+(\theta_\infty-1)^2\over
  2\nu^2},~~~\sigma=1+i\nu,~~~\nu\in{\bf R} 
$$  
If we choose the minus sign in $A$, then $iA+B=1$, and:  
$$
{y(x)= 
\left\{
1- {\nu^2+(\theta_\infty-1)^2\over \nu^2} \sin^2\left(
i {1-\sigma \over 2} \ln x +{\phi \over 2} -{\pi \over 4}
\right) ~+\delta^*(x)
\right\}^{-1} 
}
$$

$$
\phi= i\ln{4r\nu^2 \over
  i(\sigma-1)(\nu^2+(\theta_\infty-1)^2)}. 
$$

\vskip 0.3 cm 
\noindent
If we choose the plus sign in $A$, then 
$iA+B=-(\theta_\infty^\prime-1)^2/ \nu^2$, and: 
$$
y(x)= 
\left\{
 {\nu^2+(\theta_\infty-1)^2\over \nu^2} \sin^2\left(
i {1-\sigma\over 2} \ln x +{\varphi \over 2} -{\pi \over 4}
\right) ~-{(\theta_\infty-1)^2\over \nu^2}~+\delta^*(x)
\right\}^{-1} 
$$
$$
\varphi= \phi +(2k+1)\pi, ~~~k\in{\bf Z}
$$
The two ways of writing $y^\prime(x)$ give the same solution (verify
using $\sin^2=1-\cos^2$).


\section{Behaviors at $x=1,\infty$. Connection problem}
\label{conection}

 In this section is computed the behavior at $x=1$ and $x=\infty$ 
 of a solution with $p_{x1}<-2$ (i.e. $\Re \sigma_1=1$) and
 $p_{01}<-2$ (i.e. $\Re \sigma_\infty =1$) respectively. Also  
 the formulae which allow to solve the
 connection problem are computed. The results are in Proposition \ref{prop6} and
 Proposition \ref{prop7}. 

In order to understand the results, it is necessary to review the
general scheme and formulae to solve the connection
problem for $0\leq \Re \sigma_i <1$.   
In doing this, for the first time in the literature the
general 
formulas are given  (namely, the coefficients ${\bf G}_i$ in (\ref{COPPIAA}))  
 which express the monodromy data associated to a solution, 
in terms of the coefficients
of PVI and of the integration constants of the solution.  
 
\subsection{Formulae of the 
Relation between Monodromy Data and Integration Constants}

 The integration constants $\sigma$ and $r$ in (\ref{giusta}) and
 (\ref{mifido})  are functions of the
 monodromy data. These functions are computed in 
 \cite{Jimbo}. Due to a miss print  in \cite{Jimbo},  
 the correct expression  is 
 re-computed in 
 \cite{Boalch}, and the result is as follows: 
$$
2\cos\pi\sigma = p_{0x}
$$

$$
r=r(\theta_0,\theta_x,\theta_1,\theta_\infty;\sigma,p_{01},p_{x1})
$$
\be
= {(\theta_0-\theta_x+\sigma)(\theta_0+\theta_x-\sigma)(\theta_\infty+
\theta_1-\sigma)\over
 4\sigma(
 \theta_\infty+\theta_1+\sigma)} ~ {1\over {\bf F}}, 
\label{bfa}
\ee
where
$$
{\bf F}:=  {
\Gamma(1+\sigma)^2\Gamma \left({1\over 2}(\theta_0+\theta_x-\sigma)+1
\right) \Gamma\left(  
{1\over 2} ( \theta_x-\theta_0-\sigma)+1
\right)
\over 
\Gamma(1-\sigma)^2 \Gamma \left({1\over2}(\theta_0+\theta_x+\sigma)+1
\right) \Gamma\left(   
{1\over 2} ( \theta_x-\theta_0+\sigma)+1
\right)
}~\times
$$
$$
\times  {
\Gamma \left({1\over 2}(\theta_\infty+\theta_1-\sigma)+1 \right)
 \Gamma\left(  
{1\over 2} ( \theta_1-\theta_\infty-\sigma)+1
\right)
\over
\Gamma \left({1\over 2}(\theta_\infty+\theta_1+\sigma)+1 \right)
 \Gamma\left( 
{1\over 2} ( \theta_1-\theta_\infty+\sigma)+1
\right)
}~ {V\over U},
$$
and:  
$$
   U:=\left[
{i \over 2}\sin(\pi\sigma)p_{x1}- \cos(\pi\theta_x) \cos(\pi \theta_\infty) - \cos(\pi\theta_0) \cos(\pi \theta_1) \right]e^{i\pi\sigma}~+$$
$$+
{i\over 2} \sin(\pi\sigma)p_{01}+ \cos(\pi\theta_x) \cos(\pi \theta_1) + \cos(\pi\theta_\infty) \cos(\pi \theta_0)
$$
$$ 
V:=4 \sin{\pi\over 2} (\theta_0+\theta_x-\sigma) \sin{\pi\over 2} (\theta_0 - \theta_x+\sigma)~ \sin{\pi\over 2} (\theta_\infty+\theta_1-\sigma)
 \sin{\pi\over 2} (\theta_\infty-\theta_1+\sigma).
$$

\vskip 0.3 cm
\noindent
{\bf Remarks:}

{\bf 1)} Formula (\ref{bfa}) was computed with the assumption that $\theta_0,\theta_x,\theta_1,\theta_\infty,\sigma$ are not integers, and 
$\sigma\pm(\theta_0+\theta_x)$, $\sigma\pm(\theta_0-\theta_x)$, 
$\sigma\pm(\theta_1+\theta_\infty)$,
$\sigma\pm(\theta_1-\theta_\infty)$  are not even integers.   Formula (\ref{bfa}) has finite non vanishing limit when  $\sigma$ tends to  $\pm(\theta_0+\theta_x)$, $\pm(\theta_0-\theta_x)$. The corresponding solutions are (\ref{passero}) and (\ref{aquila}).

{\bf 2)} In the case $\theta_0=\theta_x=\theta_1=0$,  $r$ is computed
in \cite{DM} for the generic case, and in 
\cite{D4} {\it for all possible values of} $\theta_\infty \neq 0$ and
$\sigma\not\in (-\infty,0)\cup[1,\infty)$.

\vskip 0.3 cm
As for logarithmic solutions and Taylor expanded solutions, the parametrization of the integration constant $r$ in temrs  of monodromy data is given in  \cite{D2} \cite{D1} (and also in \cite{Jimbo} for the $\tau$ function of the logarithmic case). In this cases,  $\sigma=0$ ($p_{ij}=2$) and $\sigma=1$. 
($p_{ij}=-2$), Please, refer to these papers for the results. Here, we  concentrate on Jimbo's solutions and the sine-type oscillatory behaviors, for which $\sigma\neq 0,1$.

\vskip 0.3 cm 
 In order to solve the connection problem, also 
the inverse formulae of (\ref{bfa}) are necessary, which  give $p_{0x},p_{01},p_{x1}$ in
terms of $\sigma$, $r$ and the coefficients of PVI, namely
$\theta_0,\theta_x,\theta_1,\theta_\infty$.  To compute $p_{0x},p_{01},p_{x1}$, one has to starts from the  
monodromy matrices, which have been computed in \cite{Jimbo} for the first
time, and subsequently in 
\cite{D3}, \cite{Boalch}. Taking their traces, we obtain:  
\be
\label{COPPIAA}
\left\{\matrix{p_{0x}&=& 2\cos \pi \sigma ~~~~~~~~~~~~~~~~~~~~~~~~~~~
\cr
\cr
p_{x1}&=&{\bf G}_1 r^{-1}+{\bf G}_2+{\bf G}_3r~~~~~~~~~~~~~
\cr
\cr
p_{01}&=&{\bf G}_4r^{-1}+{\bf G}_5+{\bf G}_6 r  ~~~~~~~~~~~~~
}
\right.
\ee
where 
${\bf G}_i$ are rational functions of $\theta_\mu\pm\theta_\nu\pm\sigma$, $\cos{\pi\sigma}$,
$\cos{\pi\theta_\mu}$, $\Gamma((\theta_\mu\pm\theta_\nu\pm\sigma)/2)$,
$e^{\pm i\pi\sigma}$.
 The explicit computation
of the coefficients ${\bf G}_i$'s is very complicated, and it is not written anywhere (except for \cite{DM}
and \cite{D4},
 when $\theta_0=\theta_x=\theta_1=0$). So, here the $G_i$'s are given, for the first
 time. To do the computation,  the monodromy matrices of   \cite{D3}, page
1355-1357, formulae (A23), (A24), (A25) ($r$ appears with the name $s$) are used, with the assumption that  $\theta_0,\theta_x,\theta_1,\theta_\infty,\sigma$ are not integers. Here is the result. 
  \vskip 0.2 cm 
\noindent
Let $s(z):=\sin({\pi\over 2}z)$ and 

$$
\Xi=\Bigl(s(\theta_0+\theta_x+\sigma)s(\theta_0-\theta_x-\sigma)
+s(\theta_0-\theta_x+\sigma)s(\theta_0+\theta_x-\sigma)\Bigr)\times
$$
$$
\times
\Bigl(s(\theta_1+\theta_\infty+\sigma)s(\theta_1-\theta_\infty+\sigma)
+s(\theta_1+\theta_\infty-\sigma)s(\theta_1-\theta_\infty-\sigma)\Bigr)
$$

$$
\Xi_1=\Bigl(s(\theta_0+\theta_x+\sigma)s(\theta_0-\theta_x+\sigma)
+s(\theta_0+\theta_x-\sigma)s(\theta_0-\theta_x-\sigma)\Bigr)\times
$$
$$
\times
\Bigl(s(\theta_1+\theta_\infty+\sigma)s(\theta_1-\theta_\infty+\sigma)
+s(\theta_1+\theta_\infty-\sigma)s(\theta_1-\theta_\infty-\sigma)\Bigr)
$$

$$
\Omega=\Bigl(-s(\theta_0+\theta_x+\sigma)s(\theta_0-\theta_x-\sigma)
+s(\theta_0-\theta_x+\sigma)s(\theta_0+\theta_x-\sigma)\Bigr)\times
$$
$$
\times
\Bigl(s(\theta_1+\theta_\infty+\sigma)s(\theta_1-\theta_\infty+\sigma)
-s(\theta_1+\theta_\infty-\sigma)s(\theta_1-\theta_\infty-\sigma)\Bigr)
$$


$$
\Omega_1=\Bigl(s(\theta_0+\theta_x+\sigma)s(\theta_0-\theta_x+\sigma)
-s(\theta_0+\theta_x-\sigma)s(\theta_0-\theta_x-\sigma)\Bigr)\times
$$
$$
\times
\Bigl(s(\theta_1+\theta_\infty+\sigma)s(\theta_1-\theta_\infty+\sigma)
-s(\theta_1+\theta_\infty-\sigma)s(\theta_1-\theta_\infty-\sigma)\Bigr)
$$
Let also: 
$$
{\cal F}:=  {
\Gamma(1+\sigma)^2\Gamma \left({1\over 2}(\theta_0+\theta_x-\sigma)+1
\right) \Gamma\left(  
{1\over 2} ( \theta_x-\theta_0-\sigma)+1
\right)
\over 
\Gamma(1-\sigma)^2 \Gamma \left({1\over2}(\theta_0+\theta_x+\sigma)+1
\right) \Gamma\left(   
{1\over 2} ( \theta_x-\theta_0+\sigma)+1
\right)
}~\times
$$
$$
\times  {
\Gamma \left({1\over 2}(\theta_\infty+\theta_1-\sigma)+1 \right)
 \Gamma\left(  
{1\over 2} ( \theta_1-\theta_\infty-\sigma)+1
\right)
\over
\Gamma \left({1\over 2}(\theta_\infty+\theta_1+\sigma)+1 \right)
 \Gamma\left( 
{1\over 2} ( \theta_1-\theta_\infty+\sigma)+1
\right)
}~ {4\sigma(\theta_\infty+\theta_1+\sigma)\over (\theta_0-\theta_x+\sigma)(\theta_0+\theta_x-\sigma)(\theta_\infty+\theta_1-\sigma)}
$$
and: 
$$ 
V:=4 \sin{\pi\over 2} (\theta_0+\theta_x-\sigma) \sin{\pi\over 2} (\theta_0 - \theta_x+\sigma)~ \sin{\pi\over 2} (\theta_\infty+\theta_1-\sigma)
 \sin{\pi\over 2} (\theta_\infty-\theta_1+\sigma).
$$
$$
V_1:= V(\sigma\mapsto-\sigma)
$$
 The result is:

$$
{\bf G}_2= {2(\Omega~\cos\pi
\theta_x\cos\pi\theta_1~-~\Xi~\sin\pi\theta_x \sin\pi\theta_1)
\over
\sin^2(\pi\sigma)\sin\pi\theta_x\sin\pi\theta_1} ;
$$

$${\bf G}_5= 2\Bigl(\cos\pi\theta_1\cos\pi\theta_0+{\Xi_1\over
  \Omega_1}~\sin\pi\theta_1 \sin\pi\theta_0\Bigr);
$$
and
$$
{\bf G}_1=-{\sin\pi\theta_x \sin\pi\theta_1 \over \Omega} ~V_1~{1\over {\cal F}},~~~~~{\bf G}_3=-{\sin\pi\theta_x \sin\pi\theta_1 \over \Omega} ~V~{\cal F}
$$
$$
{\bf G}_4= - e^{i\pi\sigma}{\sin\pi\theta_0\over \sin\pi\theta_x}{\Omega\over \Omega_1} {\bf G}_1,~~~~~{\bf G}_6= - e^{-i\pi\sigma}{\sin\pi\theta_0\over \sin\pi\theta_x}{\Omega\over \Omega_1} {\bf G}_3.
$$
 Observe that the limit of the ${\bf G}_i$'s exists also for $\sigma \pm(\theta_0\pm\theta_x)\to 2n$, $n\in {\bf Z}$, though this is not always the case for the solution (\ref{bfa}) (which has anyway limit for $\sigma \pm(\theta_0\pm\theta_x)\to 0$. 
\vskip 0.2 cm 
\noindent

\subsection{Critical behaviors at $x=1,\infty$ from the behavior at $x=0$}

\vskip 0.3 cm 
 One can avoid recomputing the critical behaviors at 
$x=1,\infty$. They can be deduced from the behaviors at $x=0$  with two fractional
linear transformations, which are symmetries of PVI. 
 
The transformation $\sigma_{01}$ exchanges the values $0$ and $1$ of the
independent variable:  
\be
\sigma_{01}:~~~\theta_0^{\prime}=\theta_1,~~\theta_x^{\prime}=\theta_x,~~\theta_1^{\prime}=\theta_0,~~\theta_\infty^{\prime}=\theta_\infty;~~~y^{\prime}(x^\prime)=1-y(x),~~~x^\prime=1-x. 
\label{o-nara}
\ee
Therefore, when  $x\to 0$ then  $x^\prime\to 1$. We obtain the behavior of
$y^\prime(x^\prime)$ at $x^\prime=1$ from that of $y(x)$ at $x=0$.  
The monodromy data change  as follows \cite{D1}: 
$$
\left\{
\matrix{
p_{01}^\prime
&=&-p_{01}-p_{0x}p_{x1}+p_\infty p_x+p_1p_0 
~~~~~~~~~~~~~~~~~~~~~~~~~~~~~~~~~~~~~
\cr
\cr
p_{0x}^\prime&=&p_{x1}~~~~~~~~~~~~~~~~~~~~~~~~~~~~~~~~~~~~~~~~~~~~~~~~~~~~~~~~~~~~~~~~~~~~
\cr
\cr
p_{x1}^\prime&=&p_{0x}~~~~~~~~~~~~~~~~~~~~~~~~~~~~~~~~~~~~~~~~~~~~~~~~~~~~~~~~~~~~~~~~~~~~
}
\right.
$$
and the inverse:
\be
\label{pp}
\left\{
\matrix{
p_{01}&=&-p_{01}^\prime-p_{0x}^\prime p_{x1}^\prime+p_\infty^\prime
p_x^\prime+p_1^\prime p_0^\prime ~~~~~~~~~~~~~~~~~~~~~~~~~~~~~~~~~~~~~
\cr
\cr
p_{0x}&=&p_{1x}^\prime~~~~~~~~~~~~~~~~~~~~~~~~~~~~~~~~~~~~~~~~~~~~~~~~~~~~~~~~~~~~~~~~~~~~
\cr
\cr
p_{x1}&=&p_{0x}^\prime~~~~~~~~~~~~~~~~~~~~~~~~~~~~~~~~~~~~~~~~~~~~~~~~~~~~~~~~~~~~~~~~~~~~
}
\right.
\ee
This means
that $y^\prime$ is associated to the monodromy data with $\prime$. Namely:
$$
y^\prime(x^\prime,\Theta^\prime,P^\prime)=1-y\Bigl(x(x^\prime),\Theta(\Theta^\prime),
P(P^\prime)\Bigr)$$ 
where $\Theta$ stands for the collection of $\theta_\mu$'s, and $P$
for the collection of the $p_{ij}$'s. Formula (\ref{o-nara}) gives 
$x(x^\prime)$ and
$\Theta=\Theta(\Theta^\prime)$, while $P(P^\prime)$ is (\ref{pp}).  

\vskip 0.2 cm 
The transformation $\sigma_{x1}$ exchanges the values $x$ and $1$ of the
independent variable:  
$$
\sigma_{x1}:~~~\theta_x^\prime=\theta_1,~~\theta_1^\prime=\theta_x;~~~~~~~
\theta_0^\prime=\theta_0,~~\theta_\infty^\prime=\theta_\infty;~~~~~~~~~~
y^\prime(x^\prime)={1\over x}y(x),~~~x^\prime={1\over x}.
$$
Therefore, if $x\to 0$, $x^\prime\to \infty$ and we obtain the behavior of
$y^\prime(x^\prime)$ from that of $y(x)$.  
The monodromy data change \cite{D1} \cite{MV} as follows:
$$
\left\{
\matrix{p_{0x}^\prime &=&
 -p_{01}-p_{0x}p_{x1} +p_\infty p_x +p_0p_1 
~~~~~~~~~~~~~~~~~~~~~~~~~~~~~~~~~~~~~
\cr
\cr
p_{01}^\prime&=&p_{0x}~~~~~~~~~~~~~~~~~~~~~~~~~~~~~~~~~~~~~~~~~~~~~~~~~~~~~~~~~~~~~~~~~~~~
\cr\cr
p_{1x}^\prime&=&p_{1x}~~~~~~~~~~~~~~~~~~~~~~~~~~~~~~~~~~~~~~~~~~~~~~~~~~~~~~~~~~~~~~~~~~~~
}
\right.
$$
Namely:
\be
\label{ppp}
\left\{\matrix{p_{01} &= &-p_{0x}^\prime-p_{01}^\prime p_{x1}^\prime+p_\infty^\prime
p_1^\prime+p_0^\prime p_x^\prime ~~~~~~~~~~~~~~~~~~~~~~~~~~~~~~~~~~~~~
\cr\cr
p_{0x}&=&p_{01}^\prime~~~~~~~~~~~~~~~~~~~~~~~~~~~~~~~~~~~~~~~~~~~~~~~~~~~~~~~~~~~~~~~~~~~~
\cr\cr
p_{1x}&=&p_{1x}^\prime~~~~~~~~~~~~~~~~~~~~~~~~~~~~~~~~~~~~~~~~~~~~~~~~~~~~~~~~~~~~~~~~~~~~
}
\right.
\ee

\vskip 0.3 cm
\noindent
{\bf Remark:} The proof of (\ref{pp}) and (\ref{ppp}) see \cite{D1}. The result depends on the choice of the base of loops for the
fuchsian system associated to $y^\prime(x^\prime)$. Different choices
of loops that preserve  the ordering $1,2,3$ for $0,x^\prime,1$  correspond
to different branches of $y^\prime(x^\prime)$.
  The choice of the basis
of loops in \cite{D1}, which gives 
(\ref{pp}) and (\ref{ppp}), is actually the choice that gives the
simplest form for $p_{ij}^\prime$.  All other possible 
values of $p_{ij}^\prime$ can be obtained from (\ref{pp}) and
(\ref{ppp}) by the action of the braid group generated by (\ref{vesc}) and (\ref{vesc1}).


\subsection{Parametrization Formulae when $0<\Re\sigma_i<1$, $p_{ij}\not
\in(-\infty,-2]\cup[2,\infty)$}

\bpr
\label{prop4}
Let PVI be give, namely $\theta_0,\theta_x,\theta_1,\theta_\infty$ are given. 
Let the basis $\Gamma$ of figure \ref{figure1} be chosen, so that the monodromy data are refered to $\Gamma$.  If  $p_{0x},p_{01},p_{x1}\not
\in(-\infty,-2]\cup[2,\infty)$, there exists a transcendent  whith branches at $x=0,1,\infty$  having the following asymptotic behaviors: 
\be
\label{A}
y(x)=
\left\{\matrix{
a_{0}x^{1-\sigma_{0}}(1+\delta_0(x)),~~~x\to 0
\cr\cr
1-a_1(1-x)^{1-\sigma_{1}}(1+\delta_1(1-x)),~~~x\to 1
\cr\cr
a_\infty x^{\sigma_{\infty}}(1+\delta_\infty(x^{-1})),~~~x\to \infty
}
\right.
\ee
The branches at $x=0$, 1 and $\infty$ are identified by the following parametrization:
\be
\label{D}
 2\cos(\pi \sigma_{0})=p_{0x},~~~ 
2\cos(\pi \sigma_{1})=p_{x1},~~~ 
2\cos(\pi \sigma_{\infty})=p_{01}.
\ee
\vskip 0.2 cm   
\be
\label{B1}
a_0
=
{[\sigma_0^2-(\theta_0-\theta_x)^2][(\theta_0+\theta_x)^2-\sigma_0^2]
\over
16\sigma_0^3 r_0}
\ee
\be
\label{B2}
a_1
=
{[\sigma_1^2-(\theta_1-\theta_x)^2][(\theta_1+\theta_x)^2-\sigma_1^2]
\over
16\sigma_1^3 r_1}
\ee
\be
\label{B3}
a_\infty
=
{[\sigma_\infty^2-(\theta_0-\theta_1)^2][(\theta_0+\theta_1)^2-\sigma_{\infty}^2]
\over
16\sigma_\infty^3 r_\infty}
\ee
where:
\be
\label{C1}
 r_0=r(\theta_0,\theta_x,\theta_1,\theta_\infty;~\sigma_0,~p_{01},~~p_{x1})
\ee
\be
\label{C2}
r_1=r(\theta_1,\theta_x,\theta_0,\theta_\infty;~\sigma_1,~-p_{01}-p_{0x}p_{x1}
+p_\infty p_x+p_0p_1,
~~p_{0x})
\ee
\be
\label{C3}
r_\infty=r(\theta_0,\theta_1,\theta_x,\theta_\infty;~\sigma_\infty,~-p_{0x}-
p_{01}p_{x1}+p_\infty p_1+p_0p_x,~~p_{x1})
\ee
and $r=r(...)$ is the function (\ref{bfa}). The higher orders 
$\delta_0(x)$, $\delta_1(1-x)$, $\delta_\infty(x^{-1})$ 
 depends on their arguments as  (\ref{ddelta1}),
  with  coefficients which are certain rational functions of
$\sigma_1,r_1$ and $\sigma_\infty,r_\infty$ respectively. 
\epr

\vskip 0.3 cm
\noindent 
{\it PROOF:} The first behavior in (\ref{A})  is (\ref{giusta}). 
Second and third behaviors in (\ref{A}) are obtained
applying 
$\sigma_{01}$ and $\sigma_{x1}$ to (\ref{giusta}).
  We obtain 
$y^\prime(x^\prime)$, $x^\prime$ and then we drop the index ${}^\prime$ (prime).  

Formula (\ref{B1}) is (\ref{X}). Formula (\ref{C1}) is 
(\ref{bfa}), derived in  \cite{Jimbo} \cite{Boalch} \cite{D3} . To
obtain (\ref{B2}),  (\ref{C2}) and (\ref{B3}), (\ref{C3}), we
substitute $\theta_\nu=\theta_\nu(\theta_0^\prime, \theta_x^\prime, 
\theta_1^\prime, \theta_\infty^\prime)$ and
$p_{ij}=p_{ij}(p_{0x}^\prime,
p_{01}^\prime,p_{x1}^\prime,p_0^\prime,p_x^\prime,p_1^\prime,p_\infty^\prime)$
according to (\ref{pp}) and (\ref{ppp}) respectively. After
re-expressing the $\theta$'s and $p$'s as functions of the
$\theta^\prime$'s and $p^\prime$'s, we drop the index ${}^\prime$,
because {\it the monodromy data are the same for the given $y(x)$}. \qed

\vskip 0.2 cm 
\noindent
{\bf Remark:} The parametrization of $a_0$ and $\sigma_0$ of a branch $y(x)$ in terms of monodromy data $\theta_\mu$, $p_{ij}$ is obtained from the associated Fuchsian system for $x$ small and  $|\arg x|<\varphi\leq $, $0<\varphi\leq \pi$. The basis of loops $\Gamma$ of figure \ref{figure1} has been chosen, which produces the specific branch. The parametrizations for $a_1$ $\sigma_1$ and $a_\infty$, $\sigma_\infty$ are obtained from (\ref{pp}) (\ref{ppp}).  As pointed out in  the Remark  following (\ref{pp}) and (\ref{ppp}), they are the parametrization of a paricular branch  around $x=1$ ($|\arg (1-x)|<\varphi$) and $x=\infty$. These  are branches of the transcendent whose branch at $x=0$ has parameters $a_0$, $\sigma_0$.  The other branches are otained by the action of the braid group generated by (\ref{vesc}) and (\ref{vesc1}).


\subsection{Parametrization Formulae when  $\Re \sigma_i=0$, $p_{ij}>2$}

\bpr
\label{prop5} 
Let PVI be give, namely $\theta_0,\theta_x,\theta_1,\theta_\infty$ are given. 
Let the basis $\Gamma$ of figure \ref{figure1} be chosen, so that the monodromy data are refered to $\Gamma$.  If 
the  monodromy data are such that  $p_{0x}>2$, 
$p_{x1}>2$,  $p_{01}>1$, there exists a transcendent with branches having behavior: 
$$
y(x)= x\left\{
-2iA_0\sin^{2}\left(
i{\sigma_0\over 2}\ln x +{\phi_0\over 2}-{\pi \over 4}
\right)~+iA_0+B_0+\delta_0^{*}(x)
\right\},~~~x\to 0
$$
 
$$
y(x)= 1-(1-x)\left\{
-2iA_1\sin^{2}\left(
i{\sigma_1\over 2}\ln (1-x) +{\phi_1\over 2}-{\pi \over 4}
\right)~+iA_1+B_1+\delta_1^{*}(x)
\right\},~~~x\to 1
$$

$$
y(x)=
-2iA_\infty\sin^{2}\left(-
i{\sigma_\infty\over 2}\ln x +{\phi_\infty\over 2}-{\pi \over 4}
\right)~+iA_\infty+B_\infty+\delta_\infty^{*}(x^{-1}),~~~x\to \infty
$$
The branches are identified by the parametrization:
$$
 2\cos(\pi \sigma_{0})=p_{0x},~~~ 
2\cos(\pi \sigma_{1})=p_{x1},~~~ 
2\cos(\pi \sigma_{\infty})=p_{01}.
$$
$$
B_0= {\theta_0^2-\theta_x^2 + \sigma_0^2 \over 2\sigma_0^2},~~~
A=
 \sqrt{ {\theta_0^2\over \sigma_0^2} -B_0^2 },
~~~
\phi_0= i\ln {2r_0\over \sigma_0 A_0}
$$
$$
B_1= {\theta_1^2-\theta_x^2 + \sigma_1^2 \over 2\sigma_1^2},~~~
A=
 \sqrt{ {\theta_1^2\over \sigma_1^2} -B_1^2 },
~~~
\phi_0= i\ln {2r_1\over \sigma_1 A_1}
$$
$$
B_\infty= {\theta_0^2-\theta_1^2 + \sigma_\infty^2 \over 2\sigma_\infty^2},~~~
A=
 \sqrt{ {\theta_0^2\over \sigma_\infty^2} -B_\infty^2 },
~~~
\phi_\infty= i\ln {2r_\infty\over \sigma_\infty A_\infty}
$$
$$
  r_0,~r_1,~r_\infty \hbox{ as in (\ref{C1}), (\ref{C2}), (\ref{C3}). }
$$
$\delta^{*}_i(..)$ have the functional form of (\ref{ddelta2}).
\epr

\vskip 0.3 cm
\noindent 
{\it PROOF:} The behavior when $x\to 0$ is (\ref{mmfido}). 
$r_0$ is (\ref{bfa}). The behaviors at $x\to 1$ and
$x\to \infty$, and the functional dependence of the integration
constants on the monodromy data are proved as  for (\ref{A}), (\ref{B2}),
(\ref{B3}), (\ref{C2}), (\ref{C3}) starting from (\ref{mmfido}) and
(\ref{bfa}), via $\sigma_{01}$ and $\sigma_{x1}$. ~~~\qed


\subsection{Parametrization Formulae when $\Re \sigma_i =1$, $p_{ij}<-2$}

\bpr
\label{prop6} Let PVI be give, namely $\theta_0,\theta_x,\theta_1,\theta_\infty$ are given. 
Let the basis $\Gamma$ of figure \ref{figure1} be chosen, so that the monodromy data are refered to $\Gamma$.  
If $p_{0x}<-2$, $p_{x1}<-2$, $p_{01}<-2$, there exists a transcendent whith branches having behavior:
$$
y(x)=
\left\{
-2iA_0 \sin^2\left(i{1-\sigma_0 \over 2}\ln x +{\phi_0\over 2}-{\pi \over
  4}\right)~+iA_0+B_0~+\delta_0^{*}(x)
\right\}^{-1}
,~~~x\to 0
$$

 $$
y(x)=1-
\left\{
-2iA_1 \sin^2\left(i{1-\sigma_1 \over 2}\ln (1-x) +{\phi_1\over 2}-{\pi \over
  4}\right)~+iA_1+B_1~+\delta_1^{*}(1-x)
\right\}^{-1}
,~~~x\to 1
$$

$$
y(x)=
x\left\{
-2iA_\infty \sin^2\left(i{\sigma_\infty-1 \over 2}\ln x +{\phi_\infty\over 2}-{\pi \over
  4}\right)~+iA_\infty+B_\infty~+\delta_0^{*}(x^{-1})
\right\}^{-1}
~~~x\to \infty
$$
The branches are identified by the parametrization: 
$$
 2\cos(\pi \sigma_{0})=p_{0x},~~~ 
2\cos(\pi \sigma_{1})=p_{x1},~~~ 
2\cos(\pi \sigma_{\infty})=p_{01}.
$$
$$
\sigma_j = 1 +i \nu_j ,~~~\nu_j\in {\bf R},~~~j=0,1,\infty
$$

$$
B_0=
      {\nu_0^2+\theta_1^2-(\theta_{\infty}-1)^2\over 2
	\nu_0^2},~~~
 A_0=  
i\sqrt{
{(\theta_\infty -1)^2\over \nu_0^2}
+
B_0^2
},~~~\phi_0=i\ln{r_0\over (1-\sigma_0)A_0}
$$
$$
B_1=
      {\nu_1^2+\theta_0^2-(\theta_{\infty}-1)^2\over 2
	\nu_1^2},~~~
 A_1=  
i\sqrt{
{(\theta_\infty -1)^2\over \nu_1^2}
+
B_1^2
},~~~\phi_1=i\ln{r_1\over (1-\sigma_1)A_1}
$$
$$
B_\infty=
      {\nu_\infty^2+\theta_x^2-(\theta_{\infty}-1)^2\over 2
	\nu_\infty^2},~~~
 A_\infty=  
i\sqrt{
{(\theta_\infty -1)^2\over \nu_\infty^2}
+
B_\infty^2
},~~~\phi_\infty=i\ln{r_\infty\over (1-\sigma_\infty)A_\infty}
$$

$$
 r_0={\cal R}(\theta_0,\theta_x,\theta_1,\theta_\infty;\sigma_0,p_{01},p_{x1})
$$ 
$$
r_1={\cal R}(\theta_1,\theta_x,\theta_0,\theta_\infty;\sigma_1,-p_{01}-p_{0x}p_{x1}
+p_\infty p_x+p_0p_1,~
p_{0x})
$$
$$
r_\infty={\cal R}(\theta_0,\theta_1,\theta_x,\theta_\infty;\sigma_\infty,-p_{0x}-
p_{01}p_{x1}+p_\infty p_1+p_0p_x,~p_{x1})
$$

\be
\label{calR}
{\cal R}(\theta_0,\theta_x,\theta_1,\theta_\infty;\sigma,p_{01},p_{x1})
= {(\theta_\infty-\theta_1-\sigma)(\theta_\infty+\theta_1-2+\sigma)
(\theta_0+\theta_x+\sigma)\over
 4(1-\sigma)(
 \theta_0+\theta_x+2-\sigma)} ~ {1\over  {\bf F}^*}, 
\ee
and
$$
{\bf F}^*:=  {
\Gamma(2-\sigma)^2
\Gamma \left({1\over 2}(\theta_\infty+\theta_1+\sigma)\right) 
\Gamma\left(  {1\over 2} ( \theta_1-\theta_\infty+\sigma)+1\right)
\over 
\Gamma(\sigma)^2 
\Gamma \left({1\over2}(\theta_\infty+\theta_1-\sigma)+1\right) 
\Gamma\left( {1\over 2} ( \theta_1-\theta_\infty-\sigma)+2\right)
}~\times
$$
$$
\times  {
\Gamma \left({1\over 2}(\theta_0+\theta_x+\sigma)+1 \right)
\Gamma\left({1\over 2} ( \theta_x-\theta_0+\sigma)\right)
\over
\Gamma \left({1\over 2}(\theta_0+\theta_x-\sigma)+2 \right)
\Gamma\left( {1\over 2} ( \theta_x-\theta_0-\sigma)+1\right)
}~
 {{\cal V}\over {\cal U}},
$$
and:  
$$
  {\cal U}:=-e^{-i\pi\sigma}\left[
{i \over 2}\sin(\pi\sigma)p_{1x}+ \cos(\pi\theta_x) \cos(\pi \theta_\infty) + \cos(\pi\theta_0) \cos(\pi \theta_1) \right]~+$$
$$
-{i\over 2}
 \sin(\pi\sigma)p_{01}+ \cos(\pi\theta_x) \cos(\pi \theta_1) + \cos(\pi\theta_\infty) \cos(\pi \theta_0)
$$

$$ 
{\cal V}:=4 \sin{\pi\over 2} (\theta_0+\theta_x+\sigma) \sin{\pi\over 2} (\theta_0 - \theta_x-\sigma)~ \sin{\pi\over 2} (\theta_\infty+\theta_1+\sigma)
 \sin{\pi\over 2} (\theta_\infty-\theta_1-\sigma).
$$
 $\delta^*_0(x)$ is (\ref{ddelta3}) and $\delta^*_1(1-x)$,
$\delta^*_\infty(x^{-1})$ have the same functional dependence in $1-x$
and $x^{-1}$ respectively. 
\epr

\vskip 0.3 cm
\noindent
{\it PROOF:} The behavior when $x\to 0$ is (\ref{XXX}), which is
derived in section \ref{frazione}. In that section, (\ref{XXX}) is
obtained from (\ref{mmfido}) through the fractional linear 
transformation (\ref{infi}). 
Therefore ${\cal R}$  is obtained from $r$ in (\ref{bfa}), by
substituting $\theta=\theta(\theta^\prime)$,  $p=p(p^\prime)$
according to (\ref{infi}), and $\sigma=1-\sigma^\prime$ according to (\ref{infiinfi}) . After substitution, the index ${}^\prime$
is dropped. 

 As in the
 proof of (\ref{A}), (\ref{B2}),  (\ref{B3}), (\ref{C2}), (\ref{C3}), 
the formulae at $x=1,\infty$ are obtained applying the fractional
 linear transformations  $\sigma_{01}$ and $\sigma_{x1}$ to the formulae
 of the behavior
 at $x=0$. 
~~~\qed

\bpr
\label{prop7}
The inverse formula of (\ref{calR}), analogous to (\ref{COPPIAA}),
  is as follows. Let 
$$
{\bf G}_i={\bf
  G}_i(\theta_0,\theta_x,\theta_1,\theta_\infty,\sigma),~~~i=1,2,3,4,5,6
$$
indicate the functional dependence of
the  coefficients ${\bf G}_i$ in (\ref{COPPIAA}). Then, let us define:    
$$
 {\bf G}^*_i={\bf
 G}_i(\theta_\infty-1,\theta_1,\theta_x,\theta_0+1,1-\sigma).
$$
Then, the analogous of (\ref{COPPIAA})  for the case $\Re\sigma=1$ is:
\be
\label{COPPIAA1}
\left\{\matrix{p_{0x}&=& 2\cos \pi \sigma ~~~~~~~~~~~~~~~~~~~~~~~~~~~~~~~~~~
\cr
\cr
p_{x1}&=&{\bf G}^*_1~ r^{-1}+{\bf G}^*_2+{\bf G}^*_3~r~~~~~~~~~~~~~~~~~~
\cr
\cr
p_{01}&=&-{\bf G}^*_4~r^{-1}~-~{\bf G}^*_5~-~{\bf G}^*_6 ~r  ~~~~~~~~~~~~~
}
\right.
\ee

\epr

\noindent
{\it PROOF:} Application of (\ref{infi}). ~~~~\qed  


\subsection{Parametrization Formulae in general}

The above propositions assume that $p_{0x}, p_{x1}, p_{01}$ are of the same tipe (namely, for example, all greater than 2 or smaller than $-2$) The mixed case is of course the one to expect, for example $p_{0x}<-2$, $p_{01}\not\in(-\infty,-2]\cup[2,\infty)$ and $p_{01}>2$.    Any other combination of values of the $p_{ij}'s$ is possible. 
Therefore: 

\vskip 0.2 cm 

{\it Let PVI be give, namely $\theta_0,\theta_x,\theta_1,\theta_\infty$ are given. 
Let the basis $\Gamma$ of figure \ref{figure1} be chosen, so that the monodromy data are refered to $\Gamma$. For the given monodromy data, there exists a transcendent with branch having the  behavior of Proposition \ref{prop4}  at $x=0$    if $p_{0x} \not\in(-\infty,-2]\cup[2,\infty)$, of Proposition \ref{prop5} if $p_{0x}>2$, and of Proposition \ref{prop6} if $p_{0x}<-2$. At $x=1$  the branch has the behavior of Proposition \ref{prop4} if $p_{x1} \not\in(-\infty,-2]\cup[2,\infty)$, of Proposition \ref{prop5} if $p_{x1}>2$, and of Proposition \ref{prop6} if $p_{x1}<-2$.
At $x=\infty$ the branch has the behavior of Proposition \ref{prop4} if $p_{01} \not\in(-\infty,-2]\cup[2,\infty)$, of Proposition \ref{prop5} if $p_{01}>2$, and of Proposition \ref{prop6} if $p_{01}<-2$. The branch is identified by the corresponding parametrizations given at $x=0,1,\infty$ respectively in Propositions \ref{prop4}, \ref{prop5}, \ref{prop6}.
}


\subsection{Solution of the Connection Problem}

 Suppose that we know the behavior of $y(x)$ at  the critical points $x=0$. 
We want to write the behavior at the
 other critical points. 

\vskip 0.3 cm 
--  From the given behavior, we extract  the exponent $\sigma_0$ and
$a_0$ (or $\phi_0$).  From $a_0$ (or $\phi_0$) we compute $r_0$.  

\vskip 0.2 cm 
-- Given $\sigma_0$ and $r_0$, we can compute 
$p_{0x},p_{x1},p_{01}$  from formulae (\ref{COPPIAA}) or
(\ref{COPPIAA1}) 
(where $r=r_0$, $\sigma=\sigma_0$). 

\vskip 0.2 cm 
-- Then, we substitute $p_{0x},p_{x1},p_{01}$ in formulae of
proposition \ref{prop4}, \ref{prop5}, or \ref{prop6}  and we obtain
$a_1$ (or $\phi_1$) and 
 $\sigma_1$, $a_\infty$ (or $\phi_\infty$) and  $\sigma_\infty$.


\section{PVI associated to a Frobenius manifolds}
\label{EAUTUNNO}
 The structure of  a semi-simple Frobenius manifold of dimension 3 is
 described by a solution of a PVI equation with $\beta=\gamma=0$,
 $\delta={1\over 2}$, which means  
 $\theta_0=\theta_x=\theta_1=0$ \cite{Dub}. 
For the  solutions (\ref{giusta}) and (\ref{mmfido}), namely $0\leq
\Re \sigma<1$, $r$ reduces to: 
$$
 r=r(0,0,0,\theta_\infty,\sigma,p_{01},p_{x1})
$$ 
$$
= {\sigma~ {\cal G}^2(\sigma,\theta_\infty){\cal
  F}^2(\sigma,\theta_\infty)\over \sin^2\pi\sigma}
\Bigl[
(1+\cos\pi\theta_\infty)(1-e^{i\pi\sigma})+{i\over 2}\sin \pi
  \sigma~(p_{01}+
p_{x1}e^{i\pi\sigma})
\Bigr]
$$
where
$$
 {\cal G}(\sigma,\theta_\infty)=
 {
4^{-\sigma}\Gamma\left({1-\sigma\over 2} \right)^2
\over 
\Gamma\left(
1-{\theta_\infty\over 2}-{\sigma\over 2}
\right)
\Gamma\left(
{\theta_\infty\over 2}-{\sigma\over 2}
\right)
}
,~~~
{\cal F}(\sigma,\theta_\infty)={\cos^2\left({\pi\over 2}\sigma\right)\over
 \cos\pi\sigma-\cos\pi\theta_\infty}
$$
We remark that the above formulas hold true if $\sigma \neq 0,\pm
\theta_\infty +2m$, $m\in {\bf Z}$. In \cite{D4} we computed $r$ {\it for every case when}
$\sigma\not\in (-\infty,0)\cup[1,+\infty)$ and $\theta_\infty\neq
  0$. Please, refer to \cite{D4}, page 298-301, Theorem 2.

The connection problem is solved as in the general case

\vskip 0.5 cm
\noindent
$\diamond$ 
We now  consider a solution with $\Re \sigma_0 = 1$, namely
$p_{0x}<-2$.  
 This
special case is derived in subsection \ref{appendicite}. The critical
behavior is: 
$$
y(x)= 
\left\{
1- {\nu_0^2+(\theta_\infty-1)^2\over \nu_0^2} \sin^2\left(
 {\nu_0 \over 2} \ln x +{\phi_0 \over 2} -{\pi \over 4}
\right) ~+\delta_0^*(x)
\right\}^{-1} 
,~~~x\to 0.
$$

$$
2\cos\pi\sigma_0=p_{0x}<-2,~~~
\sigma_0= 1+ i\nu_0,
~~~
\phi_0= i\ln \left({4r_0\nu_0 \over
  \nu_0^2+(\theta_\infty-1)^2}\right)+\pi
$$
  
$$
r_0={\cal R}(0,0,0,\theta_\infty,\sigma_0,p_{01},p_{x1})=~~~~~~~~~~~~~~~~~~~~~~~~~~~~~~~~~~~~~~~~~~~~~~~~~~~~~~~~~~~
$$
$$
=
{16^{\sigma_0} 
\Gamma^2\left(1+{1\over 2}(\theta_\infty-\sigma_0)\right)
\Gamma^2\left(2-{1\over 2}(\theta_\infty+\sigma_0)\right)
\over
4 (1-\sigma_0)^3~ (\sin\pi\sigma_0)^2~ \Gamma^4\left({1-\sigma_0\over 2}\right)
}~\times
$$
$$~~~~~~~~~~~~~~~~~~~~~~~~~~~~~~~~~~~~~\times~\left[
(1+\cos\pi\theta_\infty)(1-e^{-i\pi\sigma_0})-{i\over
    2}\sin\pi\sigma_0~(p_{01}+p_{x1}e^{-i\pi\sigma_0})
\right]
$$

\vskip 0.3 cm
\noindent
If also $p_{x1}<-2$, then $y(x)$ has behavior:
$$
y(x)=1-
\left\{
1- {\nu_1^2+(\theta_\infty-1)^2\over \nu_1^2} \sin^2\left(
{\nu_1 \over 2} \ln(1- x) +{\phi_1 \over 2} -{\pi \over 4}
\right) ~+\delta_1^*(1-x)
\right\}^{-1} 
,~~~x\to 1.
$$

$$
2\cos\pi\sigma_1=p_{x1}<-2,~~~
\sigma_1= 1+ i\nu_1,
~~~
\phi_1= i\ln\left({4r_1\nu_1 \over
  \nu_1^2+(\theta_\infty-1)^2}\right)+\pi
$$

$$
r_1={\cal R}(0,0,0,\theta_\infty;\sigma_1,-p_{01}-p_{0x}p_{x1}
+4(\cos(\pi\theta_\infty)+1),
p_{0x})=
~~~~~~~~~~~~~~~~~~~~~~~~~~~~
$$
$$
=
{16^{\sigma_1} 
\Gamma^2\left(1+{1\over 2}(\theta_\infty-\sigma_1)\right)
\Gamma^2\left(2-{1\over 2}(\theta_\infty+\sigma_1)\right)
\over
4 (1-\sigma_1)^3~ (\sin\pi\sigma_1)^2~ \Gamma^4\left({1-\sigma_1\over 2}\right)
}~\times~~~~~~~~~~~~~~~~~~~~~~~~~~~~~~~~~~~~~~~~
$$
$$~~~~~~~~~~~~~~~~~~~~~~~~~~~~~~~~~~~~~\times~\left[
(1+\cos\pi\theta_\infty)(1-e^{i\pi\sigma_1})-{i\over
    2}\sin\pi\sigma_1~(p_{0x}e^{-i\pi\sigma_1}-p_{01}-p_{0x}p_{x1})
\right]
$$

\vskip 0.3 cm
\noindent
If also $p_{01}<-2$, then $y(x)$ has behavior:
$$
y(x)=x
\left\{
1- {\nu_\infty^2+(\theta_\infty-1)^2\over \nu_\infty^2} \sin^2\left(
- {\nu_\infty \over 2} \ln x +{\phi_\infty \over 2} -{\pi \over 4}
\right) ~+\delta_\infty^*(1-x)
\right\}^{-1} 
,~~~x\to \infty.
$$

$$
2\cos\pi\sigma_\infty=p_{01}<-2,~~~
\sigma_\infty= 1 + i\nu_\infty,
~~~
\phi_\infty= i\ln\left({4r_\infty\nu_\infty \over
  \nu_\infty^2+(\theta_\infty-1)^2}\right)+\pi
$$

$$
r_\infty={\cal R}(0,0,0,\theta_\infty;\sigma_\infty,-p_{0x}-
p_{01}p_{x1}+4(\cos(\pi\theta_\infty)+1),p_{x1})=
~~~~~~~~~~~~~~~~~~~~~~~~
$$
$$
=
{16^{\sigma_\infty} 
\Gamma^2\left(1+{1\over 2}(\theta_\infty-\sigma_\infty)\right)
\Gamma^2\left(2-{1\over 2}(\theta_\infty+\sigma_\infty)\right)
\over
4 (1-\sigma_\infty)^3~ (\sin\pi\sigma_\infty)^2~ 
\Gamma^4\left({1-\sigma_\infty\over 2}\right)
}~\times~~~~~~~~~~~~~~~~~~~~~~~~~~~~
$$
$$~~~~~~~~~~~~~~~~~~~~~~~~~~~~~~~~~~~~~\times~\left[
(1+\cos\pi\theta_\infty)(1-e^{i\pi\sigma_\infty})-{i\over
    2}\sin\pi\sigma_\infty~(p_{x1}e^{-i\pi\sigma_\infty}-p_{0x}-p_{01}p_{x1})
\right]
$$



\section{The Full Expansion}
\label{FULL}

The full asymptotic expansion of a solution $y(x)$ when $x\to 0$,
 for $0\leq \Re
\sigma <1$, is:  
\be
\label{fullEXP}
\sq{
  y(x)= \sum_{n=1}^\infty x^n\sum_{m=-n}^n c_{nm}
  x^{m\sigma},~~~0\leq \Re \sigma<1.
}
\ee
The above series can be rigorously obtained  from the
elliptic representation of PVI of \cite{D3}. This is explained in Appendix
II, where the series follows form  $\delta_E(x)$ and  $\delta_E^*(x)$, $\phi(x)$ in (\ref{AAAa}) and (\ref{BBBb}) respectively (where $\nu_2=1-\sigma$. A special attention must be payed for the case of $\delta_E^*(x)$, $\phi(x)$. See section \ref{oscilla}). In \cite{D3} the author proved the convergence of $y(x)$  by solving some integral equations with successive approximations, derived from the elliptic representation of PVI (PVI is written as a system of two first order equations, which are solved by their associated integral equations. The solution is constructed as a series by successive approximations.   A similar procedure was first introduced by S.Shimomura (review in \cite{IKSY})). 

On the other hand, the proof of \cite{D3} does not fix the bound $-n\leq m \leq n$  (i.e. the upper bound of $m_1$ in (\ref{AAAa}), (\ref{BBBb}) must be $2m_1+1$, but it is not determined by the procedure of \cite{D3}), and it gives no recursive procedure to compute the coefficients $c_{nm}$.  This is instead possible by the recursive computational procedure explained below, by substitution of the above series into PVI. All the coefficients are
determined recursively in terms of $\sigma$ and another parameter
$r$. The series (\ref{fullEXP}) gives the series of $\delta(x)$ and $\delta^*(x)$ in Propositions \ref{prop1} and  \ref{prop2}.

 In order to compute the coefficients $c_{nm}$, we write PVI as $
Eq=0$, 
where:  
$$
Eq:=
-{d^2y \over dx^2}+{1\over 2}\left[ 
{1\over y}+{1\over y-1}+{1\over y-x}
\right]
           \left({dy\over dx}\right)^2
-\left[
{1\over x}+{1\over x-1}+{1\over y-x}
\right]{dy \over dx}
$$
$$
+
{y(y-1)(y-x)\over x^2 (x-1)^2}
\left[
\alpha+\beta {x\over y^2} + \gamma {x-1\over (y-1)^2} +\delta
{x(x-1)\over (y-x)^2}
\right]
$$
 Let us substitute the  expansion  (\ref{fullEXP}) into $Eq=0$.  We
 observe  that: 
$$
 Eq= {\hbox {numerator} \over \hbox{ denominator} },~~~
\hbox{ denominator}= 2y(1-y)(y-x)x^2(x-1)^2
$$
The denominator is not zero for $x\neq 0,1,\infty$ and $y\neq
0,x,1$. So, the coefficients are determined by 
$$
 \hbox{numerator}=0
$$
Let $c$ denote  the $c_{nm}$'s. The explicit computation gives:
$$
\hbox{numerator}= \xi_3(x,c) x^3+\xi_4(x,c) x^4+\xi_5(x,c) x^5 + ...
= \sum_{l=3}^\infty \xi_l(x,c)x^l. 
$$
where $ \xi_l(x,c)$ depends on
$c_{l-2,l-2}$, $c_{l-2,l-3}$, ..., $c_{l-2,2-l}$ and on $c_{km}$,
$k\leq l-3$.  
The first term $\xi_3(x,c)$ is: 
$$
\xi_3(x,c)=\sum_{k=-2}^2 \xi_{3k}(c)x^{k\sigma}
$$
where:
$$
\xi_{32}(c)=\xi_{32}(c_{10},c_{11}),
~~~\xi_{3k}(c)=\xi_{3k}(c_{10},c_{11},c_{1,-1}),~~~k=1,0,-1,-2.
 $$
We choose  $c_{11}$ to be the free parameter
(integration constant, the other
being $\sigma$). The coefficients 
$$
\xi_{32}(c)=\xi_{32}(c_{10},c_{11}),~~~\xi_{31}(c)=\xi_{31}(c_{10},c_{11},c_{1,-1})
$$
are linear in $c_{10},c_{1,-1}$. Then, 
$$
\xi_{32}(c_{10},c_{11})=0 ~\hbox{ determines }~ c_{10}
$$  
Substitute $c_{10}$ into $\xi_{31}(c)$. Then: 
$$
\xi_{31}(c_{10},c_{11},c_{1,-1})=0 ~\hbox{ determines }~ c_{1,-1}
$$
 For example, if we write:  
$$ c_{11}= -{r\over \sigma},$$
where $r$ is a new free parameter, we find:
 $$c_{10}=B,~~~c_{1,-1}={\sigma A^2\over 4 r},
$$ 
where:
$$
B ={\sigma^2-2\beta-1+2\delta \over 2\sigma^2} 
,
~~~
A^2+B^2=- {2\beta\over \sigma^2}. 
$$
If now we substitute $ c_{10},  c_{1,-1}$ in  
$$\xi_{30}(c_{10},c_{11},c_{1,-1}),~~~
 \xi_{3,-1}(c_{10},c_{11},c_{1,-1}),~~~
  \xi_{3,-2}(c_{10},c_{11},c_{1,-1})
$$
we verify that they vanish.  
Namely $\xi_3(x,c)=0$.

\vskip 0.2 cm 
The next step is to solve 
$$
\xi_4(x,c)=0
$$
First, we substitute  into $\xi_4(x,c)$   the integration constant
$c_{11}$ and the coefficients $c_{10},c_{1,-1}$ obtained in the
previous step. We find: 
$$
\xi_4(x,c)= \sum_{k=-4}^4 \xi_{4k}(c) x^{k\sigma}
$$
where the $\xi_{4k}$'s are linear in $c_{2m}$. Precisely:
$$
\xi_{44}(c_{2,-2})=0 ~~\hbox{ determines } ~~ c_{2,-2}
$$
$$
\xi_{43}(c_{2,-2},c_{2,-1})=0  ~~\hbox{ determines } ~~ c_{2,-1}
$$
$$
\xi_{42}(c_{2,-2},c_{2,-1},c_{20})=0  ~~\hbox{ determines } ~~ c_{20}
$$
$$
\xi_{41}(c_{2,-2},c_{2,-1},c_{20},c_{21})=0  ~~\hbox{ determines } ~~ c_{21}
$$
$$
 \xi_{40}(c_{2,-2},c_{2,-1},c_{20},c_{21},c_{22})=0  ~~\hbox{
 determines }
 ~~ c_{22}
$$
Substituting the above solutions into $\xi_4(x,c)$, we find $\xi_4(x,c)=0$,
namely $\xi_{4,-1},\xi_{4,-2},\xi_{4,-3},\xi_{4,-4}$ 
vanish on the above solutions
$c_{2m}$.

\vskip 0.2 cm 
If we proceed with $\xi_5=0$ we find again 
$$
\xi_5(x,c)= \sum_{k=-5}^5 \xi_{5k}(c) x^{k\sigma}
$$
The coefficients  $\xi_{5k}$ are  linear in $c_{3m}$. We solve 
$$
 \xi_{5k}=0,~~~k=5,4,3,2,1,0,-1
$$
and determine uniquely
$$
 c_{33},~ c_{32},~c_{31},~c_{30},~c_{3,-1},~c_{3,-2},~c_{3,-3}
$$
respectively. $\xi(x,c)$ vanishes identically on these solutions.

This is a recursive procedure to determine $y(x)$ at all orders
$x^n$. In general, 
$$
\xi_{n}(x,c)= \sum_{k=-n}^n \xi_{nk}(c) x^{k\sigma}=0
$$
determines $c_{n-2,n-2},c_{n-2,n-3},...,c_{n-2,-n+2}$ uniquely. 
The crucial point is that, for any $n$, $\xi_{nk}(c)$,
$k=n,n-1,...,4-n$  is {\it linear} in $c_{n-2,m}$  and we have a {\it finite}
number of  terms $x^{k\sigma}$, $k=-n,...,n$. 

\vskip 0.3 cm 
\noindent
We can extract the {\bf leading term} of the expansion when $0\leq \Re
\sigma <1$, to check that the result is in accordance with
propositions \ref{prop1} and \ref{prop2}: 

\vskip 0.2 cm 
-- If $0<\Re\sigma <1$: 
$$
 y(x) = ax^{1-\sigma}(1+\delta(x)),~~~~~
\delta(x)=-1+\sum_{n=0}^\infty x^n\sum_{m=-n}^{n+2}
 \tilde{c}_{nm}x^{m\sigma}
$$
$$
\tilde{c}_{nm}={c_{n+1,m-1}\over a}, ~~~a=c_{1,-1}={\sigma A^2\over
  4r}.
$$
Note that  $\delta(x)\to 0$
as $x\to 0$. We can also write 
$$
\delta(x)=  \sum_{m_2=0}^\infty \sum_{m_1=0}^{2m_2+2}
\delta_{m_1m_2}x^{m_1\sigma}x^{m_2(1-\sigma)}
 ,~~~~m_1+m_2\geq 1,~~~
\delta_{m_1m_2}=\tilde{c}_{m_2,m_1-m_2}
$$

\vskip 0.2 cm
-- If $\Re \sigma=0$: 
$$
y(x)= x
\left\{
-2iA \sin^2 \left(
i{\sigma\over 2}\ln x +{\phi \over 2}-{\pi\over 4}
\right)+iA+B+\delta(x) 
\right\},
$$

$$
\phi=i \ln{2r\over \sigma A},
$$ 

$$
\delta(x)= 
\sum_{n=1}^\infty x^n \sum_{m=-n-1}^{n+1} b_{nm} x^{m\sigma}
=
\sum_{m_2=1}^\infty \sum_{m_1=-1}^{2m_2+1} a_{m_1m_2}
x^{m_1\sigma}x^{m_2(1-\sigma)},
$$ 

$$
b_{nm}=c_{n+1,m},~~~a_{m_1m_2}=b_{m_2,m_1-m_2}.
$$

\vskip 0.2 cm
\noindent
{\bf Remarks:} 

{\bf 1)} The computation of $y(x)$ can be done without
assumptions on $\sigma$. The only condition is $\sigma\not\in {\bf
  Z}$, to avoid vanishing denominators in $c_{nm}$. If we assume
$-1 < \Re \sigma <1$, the expansions are convergent for small $x$ and $\arg x$ bounded. If we further assume  $0\leq \Re \sigma <1$ the  leading term  is as above. 

{\bf 2)} Observe that $c_{nm}\sim r^m$. Moreover, observe that he coefficients $c_{nm}$ with negative $m$ contain the factor $(\sigma^2-(\theta_0+\theta_x)^2)(\sigma^2-(\theta_0-\theta_x)^2)$. Thus, if $\sigma\in\{ \pm(\theta_0+\theta_x),\pm(\theta_0-\theta_x)\}$, these $c_{nm}$ vanish, and we have:
\be
\label{Ayafebbre} 
y(x)= \sum_{n=1}^\infty x^n\sum_{m=0}^n c_{nm} x^{m\sigma}
=
\sum_{N=0}^\infty y_N(x)~(r~x^\sigma)^N, ~~~~~\sigma\in\{ \pm(\theta_0+\theta_x),\pm(\theta_0-\theta_x)\},
\ee
where $y_N(x)$ are Taylor expansions of the form:
 $$
\left.\matrix{ 
y_0(x)=& y_1^{0}x &+& y_2^{0}x^2 &+& y_3^{0}x^3 &+&y_4^{0}x^4 &+&...
\cr
\cr
y_1(x)= &y_1^{1}x&+& y_2^{1}x^2&+& y_3^{1}x^3& +&y_4^{1}x^4 &+&...
\cr
\cr
y_2(x)= &  & &y_2^{2}x^2&+& y_3^{3}x^3&+&y_4^{2}x^4&+&...
\cr
\cr
\vdots
\cr
\cr
y_N(x)= & & & & & y_N^{N}x^N&+& y_{N+1}^{N}x^{N+1}&+&...
}\right.
$$
 The condition that $|x^{1+\sigma}|$ is the dominant term (namely $|x^{1+\sigma}|> |x^{n+m\sigma}|$, $\forall n\geq 1, 0\leq m\leq n$) is: $-1<\Re \sigma <0$.  
 The condition that $|x|$ is the dominant term (namely $|x|> |x^{n+m\sigma}|$, $\forall n\geq 1, 0\leq m\leq n$) is: 
$\Re \sigma >0$. 
 The condition that $|x^{1+\sigma}|$ is  greater than $|x^n|$, $\forall n\geq 2 $ is: 
$\Re \sigma <1$. 
Therefore, if $-1<\Re\sigma<1$, $\sigma\neq 0$, the leading terms of (\ref{Ayafebbre}) are  (\ref{passero}) and (\ref{aquila}), namely: 

$$
y(x)={\theta_0\over \theta_0+\theta_x}~x ~\mp~{r\over\theta_0+\theta_x}
~x^{1+\sigma}+...,~~~\sigma=\pm(\theta_0+\theta_x)\neq 0,
$$
$$
y(x)= {\theta_0\over \theta_0-\theta_x}~x~\mp
~{r\over\theta_0-\theta_x}~x^{1+\sigma}+...,
~~~
\sigma=\pm(\theta_0-\theta_x)\neq 0.
$$
 The higher order terms are a convergent expansion.
 We observe that formula (\ref{bfa}) has limit when $\sigma$ tends to $\pm(\theta_0+\theta_x)$, $\pm(\theta_0-\theta_x)$, so it applies here as well. If moreover the above $\sigma$ is also purely immaginary  ($\sigma=i\nu$, $\nu\in {\bf R}\backslash \{0\}$), the above expansions become (\ref{passero1}) and (\ref{aquila1}). 

If $ \Re\sigma\leq -1$, no convergence is expected. The expansion (\ref{Ayafebbre}) is proved to be convergent for $-1<\Re\sigma<1$.  We expect (but not prove here) to be convergent also for any positive $\Re\sigma$. The inequality  $|x^K|>|x^{1+\sigma}|$ holds  for $\Re\sigma >K-1$, $K\geq 1$ integer. Therefore, from (\ref{Ayafebbre}), one deduces that PVI has for positive $\Re\sigma$  two out of the four solutions of the form: 
$$
 y(x) = {\theta_0\over \theta_0+\theta_x}~x~+\sum_{n=2}^K y_n^{0} x^n ~\mp~{r\over\theta_0+\theta_x}
~x^{1+\sigma}+...,~~~\sigma=\pm(\theta_0+\theta_x)\neq 0,
$$
$$
y(x)= {\theta_0\over \theta_0-\theta_x}~x~+\sum_{n=2}^K y_n^{0} x^n~\mp
~{r\over\theta_0-\theta_x}~x^{1+\sigma}+...,
~~~
\sigma=\pm(\theta_0-\theta_x)\neq 0.
$$
The integration constant $r$ appears in the $K+1= [|\Re(\theta_0\pm\theta_x)|]+2$ term.\footnote{
Observe also that we can choose $c_{1,-1}$ as integration constant, instead of $c_{11}$. Say that we put $c_{1,-1}=\tilde{r} \in{\bf C}$. We find that $c_{nm}\sim \tilde{r}^{-m}$. This time, the $c_{nm}$ with positive $m$ have factors  $(\sigma^2-(\theta_0+\theta_x)^2)(\sigma^2-(\theta_0-\theta_x)^2)$, again leading to: 
   $$
y(x)= \sum_{n=1}^\infty x^n\sum_{m=-n}^0 c_{nm} x^{m\sigma}= \sum_{N=0}^\infty \tilde{y}_N(x)~(\tilde{r}~x^{-\sigma})^N
 ~~~\hbox{ if }\sigma\in\{ \pm(\theta_0+\theta_x),\pm(\theta_0-\theta_x)\}
$$
This is again  (\ref{Ayafebbre}).}

\vskip 0.4 cm 
\noindent
$\diamond$ 
The asymptotic expansion for $\Re \sigma =1$ is obtained from
(\ref{fullEXP}) through 
(\ref{infi}), with  the substitution
$\sigma \mapsto 1-\sigma$ (see  section \ref{frazione}): 
$$
\sq{
y(x)^{-1}= \sum_{n=0}^\infty x^n\sum_{m=-n-1}^{n+1} d_{nm}~
x^{m(1-\sigma)},~~~
\Re\sigma=1.
}
$$
Practically, to compute the coefficients $d_{nm}$, let us  call the
above solution $y^\prime(x)$,  the exponent $\sigma^\prime$ and the
parameters $\theta_\mu^\prime$. Then,
we  compute the coefficients of $y(x)$, the image of $y^\prime(x)$ via
(\ref{infi}), with
$\sigma=1-\sigma^\prime$ ($\Re \sigma=0$). 
Let $c_{nm}$ be the coefficents of the $y(x)$ in (\ref{fullEXP}). 
Then, we have:  
$$
\Bigl(y^{\prime}(x)\Bigr)^{-1}={1\over x} y(x)= { \sum_{n=1}^\infty
  x^n \sum_{m=-n}^n c_{nm} x^{m\sigma}\over x}
= 
\sum_{n=0}^\infty x^n \sum_{m=-n-1}^{n+1}
c_{n+1,m}~x^{m\sigma}. 
$$ 
This   proves that: 
$$
 d_{nm}= c_{n+1,m}
$$
We extract the leading terms. Dropping again the index $\prime$, the
final result  when $\Re
\sigma=1$ is then in accordance with proposition \ref{prop3}:
$$
y(x)^{-1}= {r\over \sigma -1}x^{1-\sigma}+B+{(1-\sigma)A^2\over 4r}
x^{\sigma-1}
+\delta^{*}(x)
$$
$$
= -2iA\sin^2\left( 
i{1-\sigma\over 2}\ln x +{\phi\over 2}-{\pi\over 4}
\right)
+iA+B +\delta^{*}(x),~~~\Re\sigma=1, 
$$
where:
$$
B={(1-\sigma)^2-2\gamma+2\alpha\over
  2(1-\sigma)^2},
~~~
A^2+B^2={2\alpha\over (1-\sigma)^2},~~~\phi=i\ln{2r\over (1-\sigma)A}
$$
and
$$
\delta^{*}(x)=\sum_{n=1}^\infty x^n \sum_{m=-n-1}^{n+1} d_{nm}~
x^{m(1-\sigma)}
~ =
\sum_{m_1=1}^\infty \sum_{m_2=-1}^{2m_1+1}
   e_{m_1m_2}x^{m_1\sigma}x^{m_2(1-\sigma)},
$$
$e_{m_1m_2}=d_{m_1,m_2-m_1}$.

\vskip 0.5 cm 
{\bf Note:} The full expansion for the logarithmic solutions can be obtained by substituting into PVI the following:
$$
y(x)= x (A_1 + B_1 \ln x + C_1 \ln^2 x + D_1 \ln^3 x + ...)+
x^2(A_2+B_2\ln x +...)+...,~~~~~x\to 0.  
$$
We obtain:     
$$
y(x)=
\left\{
\matrix{ 
 {\theta_0\over \theta_0\pm\theta_x} x + O(x^2) ~~\hbox{ [Taylor
     expansion]},
\cr\cr 
x~\left(
{\theta_0^2-B_1^2\over \theta_0^2-\theta_x^2} + B_1\ln x + 
{\theta_x^2-\theta_0^2\over 4} \ln^2 x 
\right)
+x^2(...)  +...,
\cr\cr
x~(A_1\pm \theta_0 \ln x)+x^2(...)+...,~~~\hbox{ and }\theta_0=\pm \theta_x.
}
\right.
$$
$A_1$ and $B_1$ are parameters. 
The other expansions are ontained applying the symmetries to the above.



\subsection{Full Expansions  at $x=1,\infty$}

If the three exponents $\sigma_0$, $\sigma_1$, $\sigma_\infty$ satisfy
$$
0\leq \Re \sigma_i<1, ~~~i=0,1,\infty
$$
the full expansion for $y(x)$ at the three critical points can be
computed with the symmetries $\sigma_{01}$ and $\sigma_{x1}$ of
section \ref{conection}. 
$$
y(x)=
\left\{
\matrix{
\sum_{n=1}^\infty x^n \sum_{m=-n}^n c_{nm}^{(0)}~x^{m\sigma_0},~~~x\to
0
\cr
\cr
1-\sum_{n=1}^\infty (1-x)^n \sum_{m=-n}^n c_{nm}^{(1)}~(1-x)^{m\sigma_1},~~~x\to
1
\cr
\cr
\sum_{n=0}^\infty x^{-n} \sum_{m=-n-1}^{n+1} c_{nm}^{(\infty)}~x^{-m\sigma_\infty},~~~x\to\infty
}
\right.
$$
where:
$$
  c_{nm}^{(0)}=c_{nm}^{(0)}(\sigma_0,\theta_0,\theta_x,
\theta_1,\theta_\infty,r_0)
,$$
$$
  c_{nm}^{(1)}=c_{nm}^{(0)}(\sigma_1,\theta_1,
\theta_x,\theta_0,\theta_\infty,r_1)
,~~~ 
  c_{nm}^{(\infty)}=c_{n+1,m}^{(0)}(\sigma_\infty,\theta_0,\theta_1,\theta_x,
\theta_\infty,r_\infty)
$$
See section \ref{conection} for the notations $r_0$, $r_1$, $r_\infty$.

\vskip 0.3 cm 
As we already explained, if $\Re \sigma_0=1$, the full expansion for $x\to 0$ is:
$$
y(x)= {1
\over 
        \sum_{n=0}^\infty x^n \sum_{m=-n-1}^{n+1}
	d_{nm}^{(0)}~x^{m(1-\sigma_0)}},~~~x\to 0,
$$
$$
  	d_{nm}^{(0)}=
  	d_{nm}^{(0)}(\sigma_0,\theta_0,\theta_x,\theta_1,\theta_\infty,r_0).
$$
If also $\Re \sigma_1=1$,  the full expansion for $x\to 1$ is:
$$
y(x)= 1-{1
\over 
           \sum_{n=0}^\infty (1-x)^n \sum_{m=-n-1}^{n+1} d_{nm}^{(1)}
	   ~ (1-x)^{m(1-\sigma_1)}},~~~x\to 1,
$$
$$
d_{nm}^{(1)}= d_{nm}^{(0)}(\sigma_1,\theta_1,\theta_x,\theta_0,\theta_\infty,r_1).
$$
If also $\Re \sigma_\infty=1$, the full expansion for $x\to \infty$ is:
$$
y(x)= {x
\over 
 \sum_{n=0}^\infty x^{-n} \sum_{m=-n-1}^{n+1}
	d_{nm}^{(\infty)}~x^{-m(1-\sigma_\infty)}},~~~x\to \infty,
$$
$$
  	d_{nm}^{(\infty)}=
	d_{nm}^{(0)}(\sigma_\infty,\theta_0,\theta_1,\theta_x,\theta_\infty,
r_\infty). 
$$


\section{Appendix I:  Derivation of the critical behavior when $0\leq
  \Re \sigma<1$}
\label{SatoAndAll}

\subsection{Critical Behavior of the Solution of the Schlesinger Equations}
The critical behavior follows from the Lemma 2.4.8 at page 262 of
\cite{SMW}, applied to the Schlesinger equations of the Fuchsian
system of PVI.  
Let $\hat{A}_0$, $\hat{A}_x$, $\hat{A_1}$ be independent of $x$ and
satisfy the following conditions: 
$$
\hbox{ Eigenvalues } \hat{A}_j = {\theta_j\over 2},~-{\theta_j\over
  2},~~~j=0,x,1;~~~
 \hat{A}_0+\hat{A}_x+\hat{A_1}=  -{\theta_\infty\over
 2} ~\sigma_3,
$$
We also observe that Tr$(\hat{A_0}+\hat{A_x})=0$, so the eigenvalues
have opposite sign. Let them be:    
$$
 {\sigma\over 2},-{\sigma\over 2}~:=\hbox{ eigenvalues of }\Lambda:=\hat{A}_0+\hat{A}_x.
$$

\vskip 0.2 cm 
$\diamond$ Computation of $\hat{A}_1$ and $\Lambda$. 
Suppose that $\theta_\infty\neq 0$. Let  $r_1\in{\bf C}$, $r_1\neq
0$. The condition of given eigenvalues and  the relation
$\Lambda+\hat{A}_1=-{\theta_\infty\over 2}\sigma_3$ immediately imply:
\be
\hat{A_1}= 
\pmatrix{ {\sigma^2-\theta_\infty^2-\theta_1^2\over 4 \theta_\infty} & 
-r_1
\cr
{[\sigma^2-(\theta_1-\theta_\infty)^2][\sigma^2-(\theta_1+\theta_\infty)^2]\over
  16 \theta_\infty^2 }~{1\over r_1} 
&
- 
 {\sigma^2-\theta_\infty^2-\theta_1^2\over 4 \theta_\infty} 
},
\label{hatA1}
\ee
and
\be
\Lambda=\hat{A_0}+\hat{A_x}=
\pmatrix{
{\theta_1^2-\sigma^2-\theta_\infty^2\over 4 \theta_\infty}
&
r_1
\cr
-{[\sigma^2-(\theta_1-\theta_\infty)^2][\sigma^2-(\theta_1+\theta_\infty)^2]\over
  16 \theta_\infty^2 }~{1\over r_1}
& 
-{\theta_1^2-\sigma^2-\theta_\infty^2\over 4 \theta_\infty}
}.
\label{hatA0Ax}
\ee

\vskip 0.2 cm
$\diamond$ Computation of $\hat{A}_0$ and $\hat{A}_x$. 
 For our purposes it is enough to consider the case when $\sigma\neq
 0$, so that $\Lambda$ is diagonalizable (for $\sigma=0$ see
 \cite{D1}).  Let $G_0$ be the diagonalizing matrix:
$$
 G_0^{-1} \Lambda G_0= {\sigma\over 2}~\sigma_3,~~~~~
G_0=
\pmatrix{1 & 1 
             \cr
       {(\theta_\infty+\sigma)^2 - \theta_1^2 \over 4 \theta_\infty r_1} 
&
  {(\theta_\infty-\sigma)^2-\theta_1^2\over 4 \theta_\infty r_1}
}.
$$
Let us denote:
$$
\hat{\hat{A_i}}=G_0^{-1} \hat{A}_i G_0,~~~~~i=0,x.
$$
Let $r\in{\bf C}$, $r\neq 0$. 
If $\sigma\neq 0$, we have:  
\be
\hat{\hat{A_0}}
=
\pmatrix{{\theta_0^2-\theta_x^2+\sigma^2\over 4\sigma}
&
r
\cr
-{[\sigma^2-(\theta_0-\theta_x)^2][\sigma^2-(\theta_0+\theta_x)^2]
\over 16 \sigma^2 }~{1\over r}
&
-{\theta_0^2-\theta_x^2+\sigma^2\over 4\sigma}
},
\label{hathatA0}
\ee
\be
\hat{\hat{A_x}}= \pmatrix{ {\sigma^2+\theta_x^2-\theta_0^2\over 4 \sigma} 
&
-r
\cr
{[\sigma^2-(\theta_0-\theta_x)^2][\sigma^2-(\theta_0+\theta_x)^2]\over
  16 \sigma^2 }~{1\over r}
&
- {\sigma^2+\theta_x^2-\theta_0^2\over 4 \sigma} 
}.
\label{hathatAx}
\ee
\label{elicopter0}

The  lemma 2.4.8 at page 262 of
\cite{SMW},  becomes the theorem at page 1145-1146 of \cite{Jimbo}, namely:

\ble
\label{piove1}
Suppose that  $|\Re \sigma|<1$. Choose two positive numbers $\sigma_1$
and $K$ such that:
$$
|\Re \sigma|<\sigma_1<1,~~~~~~~||\hat{A_i}||<K,~~~i=0,x,1.
$$
Then, for every $\varphi>0$ there exists $\epsilon>0$ such that the
Schlesinger equations have a unique solution $A_0(x),A_x(x),A_1(x)$
holomorphic in the sector $\{x ~|~~
0<|x|<\epsilon,~~|\hbox{arg}~x|<\varphi\}$, and satisfying the
asymptotic conditions:
$$
||A_1-\hat{A_1}||<K|x|^{1-\sigma_1},~~~~~
||x^{-\Lambda}(A_1-\hat{A_1})x^{\Lambda}||<K^2|x|^{1-\sigma_1}
$$
$$
||x^{-\Lambda}A_0 x^{\Lambda}-\hat{A_0}||<K|x|^{1-\sigma_1},
~~~~~~~||x^{-\Lambda}A_x
x^{\Lambda}-\hat{A_x}||<K|x|^{1-\sigma_1}
$$
\ele

\ble
\label{piove2}
The asymptotic behavior of $A_1$ is:
$$
  A_1(x)=\hat{A_1}+\Delta_1(x),~~~~~~~\Delta_1(x)=O(x^{1-\sigma_1}),~~~x^{-\Lambda}\Delta_1(x)x^{\Lambda}= O(x^{1-\sigma_1}).
$$

The asymptotic behaviors of $A_0$ and $A_x$ are:
$$
  A_j(x)~=
  x^{\Lambda}\hat{A_j}x^{-\Lambda}~+\Delta_j(x)~=G_0 \Bigl[x^{{\sigma\over
  2}\sigma_3} \hat{\hat{A_j}} x^{-{\sigma\over 2}\sigma_3}\Bigr] G_0^{-1}~+\Delta_j(x),
$$

$$
\Delta_j(x)=O(x^{1-\sigma_1-|\Re\sigma|}),~~~j=0,x
$$
where
$$
x^{{\sigma\over
  2}\sigma_3} \hat{\hat{A_0}} x^{-{\sigma\over 2}\sigma_3}
=
\pmatrix{{\theta_0^2-\theta_x^2+\sigma^2\over 4\sigma}
&
rx^{\sigma}
\cr
-{[\sigma^2-(\theta_0-\theta_x)^2][\sigma^2-(\theta_0+\theta_x)^2]
\over 16 \sigma^2 }~{1\over r}~x^{-\sigma}
&
-{\theta_0^2-\theta_x^2+\sigma^2\over 4\sigma}
},
$$

$$
x^{{\sigma\over
  2}\sigma_3} \hat{\hat{A_x}} x^{-{\sigma\over 2}\sigma_3}= \pmatrix{ {\sigma^2+\theta_x^2-\theta_0^2\over 4 \sigma} 
&
-rx^{\sigma}
\cr
{[\sigma^2-(\theta_0-\theta_x)^2][\sigma^2-(\theta_0+\theta_x)^2]\over
  16 \sigma^2 }~{1\over r}~x^{-\sigma}
&
- {\sigma^2+\theta_x^2-\theta_0^2\over 4 \sigma} 
}.
$$
\ele

\noindent
{\it Proof:} The behavior of $A_1$ is immediately obtained from 
 lemma \ref{piove1}.  The behaviors of $A_0$, $A_x$ follow from  lemma \ref{piove1}: 
$$
 A_j= x^{\Lambda} \hat{A_j}x^{-\Lambda}+  x^{\Lambda}\tilde{\Delta_j}
x^{-\Lambda},~~~\tilde{\Delta_j}(x)=O(x^{1-\sigma_1}),~~~~~j=0,x.
$$
Observe that:
$$
x^{\Lambda}\tilde{\Delta_j}
x^{-\Lambda}= G_0~ x^{{\sigma\over
    2}\sigma_3}~(G_0^{-1}\tilde{\Delta_j}G_0)~x^{-{\sigma\over
    2}\sigma_3}~G_0^{-1}.
$$
Since $G_0$ is constant, $ x^{{\sigma\over
    2}\sigma_3}~(G_0^{-1}\tilde{\Delta_j}G_0)~x^{-{\sigma\over
    2}\sigma_3}$ has form: 
$$
x^{{\sigma\over
    2}\sigma_3}\pmatrix{m_{11}& m_{12}\cr m_{21} & m_{22}}x^{-{\sigma\over
    2}\sigma_3}= \pmatrix{m_{11}& m_{12} x^{\sigma}\cr m_{21}x^{-\sigma} & m_{22}}.
$$
the results follows, with $\Delta_j= x^{\Lambda}\tilde{\Delta_j}
x^{-\Lambda}$ .
~~~~~~~~~~~~~~~~~~~~~~~~~~~~~~~~~~~~~~~~~~\qed

\subsubsection{ Critical Behavior of $y(x)$}
\label{maldipiedi}

As it is known, the Schlesinger equations can be written in Hamiltonian form and reduce
 to PVI, being the transcendent $y(x)$ solution of
 $A(y(x),x)_{1,2}=0$. Namely:
$$
y(x)= 
{x~(A_0)_{12} \over x~\left[
(A_0)_{12}+(A_1)_{12}
\right]- (A_1)_{12}},
$$
Lemma \ref{piove2} implies: 
$$
(A_0)_{12}= r_1\left\{
{[\sigma^2-(\theta_0-\theta_x)^2][(\theta_0+\theta_x)^2-\sigma^2]
\over
16\sigma^3 r}
~x^{-\sigma}+{\theta_0^2-\theta_x^2+\sigma^2\over 2\sigma^2}
-{r\over \sigma}~x^\sigma
\right\}~+\delta_0(x),
$$
$$
 (A_{1})_{12}= -r_1 +
 \delta_1(x),~~~\delta_0(x)=O(x^{1-\sigma_1-|\Re\sigma|}),
~~~\delta_1(x)=O(x^{1-\sigma_1})
$$
For brevity, let us write 
$(A_0)_{12}=
ax^{-\sigma}+bx^{\sigma}+c+\delta_0(x)$.   Thus:
$$
y(x)= {x(ax^{-\sigma}+bx^{\sigma}+c+\delta_0(x))
\over
x(ax^{-\sigma}+bx^{\sigma}+c-r_1+\delta_0(x)+\delta_1(x))+r_1-\delta_1(x)}.
$$

\vskip 0.2 cm 
\noindent
Observe that we can restrict to $0\leq \Re\sigma<1$, being  the
negative sigma case symmetrical. 

\vskip 0.3 cm
\noindent
$\diamond$ {\bf Case $0<\Re\sigma<1$:} When $x\to 0$, the  term $x^{-\sigma}$ is  dominant over
$\delta_0(x)$ and $\delta_1(x)$. But constant terms  and $x^\sigma$ may be
of higher order than $\delta_1(x)$ and $\delta_0(x)$. Thus:
$$
y(x)=
{x\Bigl(ax^{-\sigma}+O(x^{1-\sigma_1-\Re\sigma})+O(1)+O(x^{\Re\sigma})\Bigr)
\over r_1+O(x^{1-\sigma_1})~+
x\Bigl(ax^{-\sigma}
+O(x^{1-\sigma_1-\Re\sigma})+O(1)+O(x^{\Re\sigma})\Bigr)
}
$$
$$
=
{a x^{1-\sigma}\Bigl(
1+O(\hbox{max}\{x^{1-\sigma_1},x^{\Re\sigma}\})
\Bigr)
\over 
r_1\Bigl(
1 
+
O(\hbox{max}\{x^{1-\sigma_1},x^{1-\Re\sigma}\})
\Bigr)
}
= {a\over r_1} x^{1-\sigma} ~\Bigl(
1
+O(\hbox{max}\{x^{1-\sigma_1},x^{\Re\sigma},x^{1-\Re\sigma}\})
\Bigr)
$$
Restoring the value of $a$, we find the following critical
behavior when $0<\Re\sigma<\sigma_1<1$:
\be
\label{J1}
y(x)= 
{[\sigma^2-(\theta_0-\theta_x)^2][(\theta_0+\theta_x)^2-\sigma^2]
\over
16\sigma^3 r} x^{1-\sigma} ~\Bigl(
1
+O(\hbox{max}\{x^{1-\sigma_1},x^{\Re\sigma}\})
\Bigr).
\ee

\vskip 0.3 cm 
\noindent
$\diamond$ {\bf Case $\Re\sigma=0$, $\sigma\neq 0$:} In this case
$\delta_0(x)$ and $\delta_1(x)$ are $O(x^{1-\sigma_1})$, for any
$0<\sigma_1<1$. We can choose $\sigma_1$  as small as we like. Also note that
$x^{\pm\sigma}=O(1)$, namely it is bounded for $x\to 0$ and does not
vanish. Thus:  
$$
y(x)= { x(ax^{-\sigma}+bx^\sigma+c) ~+x\delta_0(x)
   \over
r_1\left[
1-{\delta_1(x)\over r_1} +x\Bigl(
O(1)+O(x^{1-\sigma_1})
\Bigr)
\right]}
= {x\over r_1}\Bigl(
ax^{-\sigma}+bx^{\sigma}+c +r_1\delta_0(x)
\Bigr)~(1+O(x^{1-\sigma_1}))
$$
Now, if we substitute $a,b,c$ and write $x^\sigma=\exp\{\sigma\ln
x\}$, we obtain: 
 \be
\label{J2}
 y(x) = x ~\left\{ 
iA\sin\left(i\sigma \ln x +  i\ln {2r\over \sigma A}\right) +  {\theta_0^2-\theta_x^2 + \sigma^2 \over 2\sigma^2} + \delta(x)
\right\}
\bigl( 
1+ \hat{\delta}(x)
\bigr), 
\ee
where
$$
A={ \sqrt{
\bigr[\sigma^2-(\theta_0+\theta_x)^2
\bigl]
\bigr[
(\theta_0-\theta_x)^2-\sigma^2
\bigl]
}
\over 2\sigma^2
},~~~~~
\delta(x),~ \hat{\delta}(x),~\delta^*(x) ~= O(x^{1-\sigma_1}).
$$ 
\qed

\vskip 0.2 cm 
Note that if $\sigma\in\{\pm(\theta_0+\theta_x),\pm(\theta_0-\theta_x)\}$, $A$ is zero, and the coefficient of $x^{-\sigma}$ in $(A_0)_{12}$ becomes zero. Newertheless, $y(x)$ is well defined, staring with power $x$ and $x^{1+\sigma}$ ($-1<\Re\sigma<1$). It is given by (\ref{passero}), (\ref{aquila}).  
If moreover $\sigma=i\nu$, $\nu\in {\bf R}\backslash \{0\}$, $y(x)$ becomes (\ref{passero1}), (\ref{aquila1}).

\vskip 0.2 cm 
 The leading term extracted in (\ref{J1}) holds for $0\leq \Re\sigma<1$. If instead we choose $-1 < \Re\sigma \leq 0$, we would extract the term $x^{1+\sigma}$. Suppose then that, for $\sigma$ and $\tilde{\sigma}$, with  $0\leq \Re\sigma<1$ and $-1 < \Re\tilde{\sigma} \leq 0$ respectively, we have the two solutions of a given PVI:   $y(x)\sim a x^{1-\sigma}$ and  $\tilde{y}(x)\sim \tilde{a} x^{1+\tilde{\sigma}}$. Clearly: 
$$
 a= 
{[\sigma^2-(\theta_0-\theta_x)^2][(\theta_0+\theta_x)^2-\sigma^2]
\over
16\sigma^3 r},~~~~~\tilde{a}= -{\tilde{r}\over \tilde{\sigma}}
$$
If $y$ and $\tilde{y}$ are the same branch corresponding to the same monodromy data, then  Tr$(M_0M_x)=2\cos(\pi\sigma)\equiv 2\cos(\pi\tilde{\sigma})$, namely $\tilde{\sigma}=-\sigma$, and  me must have $a=\tilde{a}$. Namely: 
\be
\label{referee}
\tilde{r} = {\sigma^2 A^2\over 4~r}
\ee
We remark that $r$, given in (\ref{bfa}), does not vanish for the values of $\sigma \in \{\pm(\theta_0+\theta_x),\pm(\theta_0-\theta_x)\}$. On the other hand, $A$ is vanishes, and so $\tilde{r}$ and the first term $ax^{1-\sigma}\equiv \tilde{a}x^{1+\tilde{\sigma}}$. This is nothing but the fact that the expansion is in this case is (\ref{passero}), (\ref{aquila}), and $\tilde{r}$ is not a good integration constant in this case (see also Remark {\bf 2)} in section \ref
{FULL}). 

\vskip 0.2 cm 
We estabilish the invariance of (\ref{mifido}) when $\sigma\mapsto -\sigma$. Observe that, for purely imaginary $\sigma$ and $\tilde{\sigma}$, we have $y(x)= x\{iA\sin(i\sigma\ln x +\phi)+B+\delta^*(x)\}$, $\tilde{y}(x)= x\{i\tilde{A}\sin(i\tilde{\sigma}\ln x +\tilde{\phi})+\tilde{B}+\delta^*(x)\}$, where $\phi=i\ln (2r/\sigma A)$ $\tilde{\phi}=i\ln (2\tilde{r}/\tilde{\sigma} \tilde{A})$. 
Again, If $y$ and $\tilde{y}$ are the same branch corresponding to the same monodromy data, we have  $\tilde{y}=y$, with $\tilde{\sigma}=-\sigma$. Clearly, $\tilde{A}=A$, $\tilde{B}=B$ and the relation  (\ref{referee}) implies $\tilde{\phi}= -\phi+(2k+1)\pi$, $k\in{\bf Z}$. This means that $\sigma \mapsto -\sigma$ leaves (\ref{J2}) (namely (\ref{mifido})) invariant. 

We estabilish the invariance of (\ref{LABEL}) when $\sigma\mapsto 2-\sigma$. This is done as above, this time observing that the role of $\sigma$ is played by $1-\sigma$, and $\tilde{r}=(\sigma-1)^2A^2/4r$, where $A$ is (\ref{AA}) (recall the construction of $y(x)$ by a symmetry transformation in Section \ref{frazione}, and recall that $\phi = i \ln (2r/(1-\sigma)A)$).  This implies that $\sigma \mapsto 2-\sigma$ induces $\tilde{\phi} =  -\phi +(2k+1)\pi$. 



\section{Appendix II: Elliptic Representation}\label{elliptic}

In this paper, all the critical behaviors are revised for any $\sigma$
such that  $0\leq \Re
\sigma \leq 1$ 
$\sigma\neq 0,1$.  In \cite{D3} all the critical 
behaviors  for any $0\leq \Re \sigma \leq 1$, 
 $\sigma\neq 0,1$, are also obtained  using the elliptic representation of
 PVI
\footnote{Actually, for  any $\sigma\not \in (-\infty,0]\cup
	 [1,\infty)$.  But  $\Re \sigma<0$ or $>1$ is equivalent to
	 $0\leq \Re \sigma \leq 1$.}
.
  If $0<\Re \sigma<1$, the behavior (\ref{giusta}) is exactly the
behavior 
(\ref{AAAa}) computed in \cite{D3}. 

But when $\Re\sigma=0,1$, the
critical behaviors  
of $y(x)$
 obtained in  \cite{D3}, namely (\ref{BBBb}) and (\ref{CCCc})
 below, are apparently  different (\ref{mmfido}) and  (\ref{XXX}). 
 Now,  (\ref{mmfido}) must
	 coincide with   
(\ref{BBBb}), and (\ref{XXX}) with (\ref{CCCc}). They are just written
	 in a different way. This coincidence allows  to prove the convergence
	 of the series of $\delta(x)$ and $\delta^*(x)$. 

\vskip 0.3 cm 

Before showing this coincidence, let us review the
	 elliptic representation of a Painlev\'e VI function. This 
is: 
\be
\label{ELIC}
  y(x) =\wp\left( \nu_1 \omega_1(x)
  +\nu_2\omega_2(x)+v(x);\omega_1,\omega^{}_2 \right)+{1+x\over
    3},~~~\nu_1,\nu_2\in{\bf C},
\ee
where $\omega_1$, $\omega_2$ are the half-periods. 
$\omega_1$ is the hypergeometric function:  
\be
\label{FffF}
\omega_1(x)={\pi \over 2} F\left({1\over 2},{1\over 2},1;x\right)
\ee
and
$$
\omega_2(x)= -{i\over 2}[F\left({1\over 2},{1\over 2},
1;x\right)\ln(x)+F_1(x)],~~~~|\arg x|<\pi 
$$
\be
\label{FffF1}
F_1(x):=  
\sum_{n=0}^{\infty}{ \left[\left({1 \over 2}\right)_n\right]^2
  \over (n!)^2 } 2\left[ \psi(n+{1\over 2}) - \psi(n+1)\right]
x^n,
\ee
$$\psi(z) = 
{d \over dz}\ln \Gamma(z),~~~ \psi\left({1\over 2}\right) = -\gamma -2 \ln
2,~~~ \psi(1)=-\gamma,~~~\psi(a+n)=\psi(a)+\sum_{l=0}^{n-1} {1\over a+l}.
$$
The function $v(x)$ solves a non linear equation equivalent to PVI,
and in \cite{D3} it is proved that it has a convergent expansion. Namely,
for any complex  $\nu_1$ and $\nu_2$ , such that $\nu_2\not \in
(-\infty,0]\cup\{1\}\cup[2,+\infty)$, there exists a sufficiently
    small $\epsilon<1$ and a solution $v(x)$ such that:  
$$
v(x) = \sum_{n\geq 1}a_n x^n+
\sum_{n\geq 0,m\geq 1} b_{nm} x^n\Bigl[e^{-i\pi
    \nu_1}x^{1-\nu_2}\Bigr]^m 
+
\sum_{n\geq 0,m\geq 1} c_{nm} x^n\Bigl[e^{i\pi
    \nu_1}x^{\nu_2}\Bigr]^m
 $$
$$
=\sum_{m_1m_2}v_{m_1m_2}x^{m_1(1-\nu_2)+m_2\nu_2},~~~m_1+m_2\geq
1,~~~m_1,~m_2\geq 0.
$$
$a_n,b_{nm},c_{nm}$  
are certain rational functions of $\alpha,\beta,\gamma,\delta,\nu_2$ . The series is proved to converge (see \cite{D3})  and defines an
holomorphic function of $x,x^{\nu_2},x^{1-\nu_2}$ in the domain:
$$
{\cal D}= \{ x\in\tilde{{\bf C}\backslash\{0\}}~|~|x|<\epsilon, ~|e^{i\pi
    \nu_1}x^{\nu_2}|<\epsilon,~|e^{-i\pi
    \nu_1}x^{1-\nu_2}|<\epsilon
\},
$$
$$
{\cal D}= \{  x\in\tilde{{\bf C}\backslash\{0\}}~|~|x|<\epsilon \} \hbox{ if } \Im
\nu_2=0. 
$$
 The critical behavior will be 
determined by the exponent $\nu_2$, which is identified with $\sigma$ 
in the following way: 
$$
  \nu_2=1-\sigma\hbox{ if } \nu_2\not \in(1,2)
$$
$$
   \nu_2=1+\sigma \hbox{ if } \nu_2 \in(1,2)
$$
where $(a,b)$ is the notation for an open interval.  

The asymptotic behavior of $y(x)$ is obtained from the Fourier
expansion of the $\wp$-function. Let the 
``modular parameter''  be:
$$ 
\tau(x)= {\omega_2(x)\over \omega_1(x)} = {1\over \pi }
\left(
\arg x - i\ln {|x|\over 16}
\right)-{i\over \pi} \left(
{F_1\over F}+\ln16
\right)
$$
Note that $F_1/F+\ln 16=O(x)$. 
The elliptic function can be expanded, when $x\to 0$, as a convergent
Fourier series, under the condition (satisfied in ${\cal D}$) that: 
$$
\Im \tau \geq \left|\Im\left({ \nu_1\omega_1+\nu_2\omega_2+v \over 2\omega_1}\right)
\right|.
$$
 The expansion is: 
 $$
  y(x) = \wp(\nu_1\omega_1+\nu_2\omega_2+v;\omega_1,\omega_2)+{1+x\over 3}
=
$$
\be
\label{uccellogiallo}
\left({\pi \over 2 \omega_1}\right)^2
\left\{ 
-{1\over 3} +\sin^{-2}\left( 
{f \over 2} 
\right)
+8  \sum_{n\geq 1} {n e^{2 i \pi  n \tau} \over 1 - e^{2 i \pi  n \tau} }
\left[ 
1-\cos
\left(
 n 
f 
\right)
\right]
\right\}
              +{1+x\over 3}
\ee
 where $f:= 
  \nu_1 +\nu_2\tau +{v\over \omega_1}
$. 
Note that in ${\cal D}$,  
$\left| e^{i f(x)}\right|<1$  and $\sin\left({\pi \over 2}
  f\right)\neq 0$. Namely, 
the denominator in the expansion does not vanish in ${\cal D}$. 

\vskip 0.3 cm

Now let us consider the case  $0\leq \Re \sigma \leq 1$, $\sigma\neq 0,1$,
 namely $ 0\leq \Re \nu_2 \leq 1$, $\nu_2\neq 0,1$. In this case,
 $$
{\cal D}=\{ x~ | ~0<|x|<\epsilon\}.
$$
 The other cases (namely, $\Re
 \sigma<0$, $\Re\sigma>1$, $\sigma\not\in (-\infty,0]\cup[1,\infty)$) are
 equivalent to the above, as it is proved in \cite{D3}.


\subsection{ Case $0< \Re \nu_2 \leq 1$, namely $0\leq \Re \sigma<1$}

We expand (\ref{uccellogiallo}) when $x\to 0$, keeping dominant terms:
$$ 
y(x) = {x\over 2} - 4 e^{i\pi \nu_1} \left({x\over 16}\right)^{\nu_2} 
e^{i{\pi v(x)\over \omega_1}}-4 e^{-i\pi \nu_1} \left({x\over 16}\right)^{2-\nu_2} 
e^{-i{\pi v(x)\over \omega_1}} +~~~~~~~~~~~~~
$$
$$~~~~~~~~~~~~~~~~
  +O(\max\{x~x^{\nu_2},~x~x^{2-\nu_2},~x^{2\nu_2},~x^2,
~x^{4-2\nu_2}\})
$$

\vskip 0.2 cm 
\noindent
$\diamond$ {\bf Case $0< \Re \nu_2 < 1$, namely $0< \Re \sigma<1$:} In
this case $v(x)\to 0$ for $x\to 0$, $
 e^{i\pi {v(x)\over \omega_1(x)}}= 1+
 O(x)+O(x^{\nu_2})+O(x^{1-\nu_2})$. 
 From the expansion of $v(x)$ and (\ref{uccellogiallo}) we compute:
\be
\label{AAAa}
{
 y(x)= -4e^{i\pi\nu_1} \left({x\over
 16}\right)^{\nu_2}(1+\delta_E(x))
,~~~\nu_2=1-\sigma}.
\ee
$$
\delta_E(x)= 
\sum_{m_1\geq 0,m_2\geq 0,m_1+m_2\geq 1}\delta_{m_1m_2}x^{m_1(1-\nu_2)+m_2\nu_2}
=O(\max\{x^{\nu_2},x^{1-\nu_2}\}).
$$
$\delta_{m_1m_2}\in{\bf C}$.
This behavior  coincides with (\ref{giusta}). The series
$\delta_E(x)$ converges in ${\cal D}$ and 
 coincides with (\ref{ddelta1}). This
proves the convergence of   (\ref{ddelta1}). 

\vskip 0.2 cm 

\vskip 0.2 cm 
\noindent
{\it Remark:} 
For $1<\nu_2<2$, we obtain:
$$
y(x)= -4e^{-i\pi\nu_1}\left({x\over
  16}\right)^{2-\nu_2}(1+O(\max\{x^{2-\nu_2},x^{\nu_2-1}\}) ),~~~\nu_2=1+\sigma
$$

\vskip 0.2 cm 
\noindent
$\diamond$ {\bf Case $ \Re\nu_2 = 1$, i.e. $\Re\sigma=0$} Now $v(x)\not\to 0$, namely:  
$$ 
v(x)= \phi(x) + O(x), ~~~~~
\phi(x):= 
\sum_{m\geq 1} 
b_{0m} \left[ e^{-i\pi \nu_1} x^{1-\nu_2}\right]^m
\not \to 0 \hbox{ as } x\to 0
$$
and $e^{i\pi {v\over \omega_1}}= e^{2i\phi}(1+O(x))$. The series of
$\phi(x)$ converges in ${\cal D}$.  
The dominant terms in the Fourier expansion are (note that $x$, $x^{\nu_2}$ and 
$x^{2-\nu_2}$ are of the same order):
$$ 
y(x) = {x\over 2} 
   - 4 e^{i\pi \nu_1} \left({x\over 16}\right)^{\nu_2} 
e^{i{\pi v(x)\over \omega_1}}
   -4 e^{-i\pi \nu_1} \left({x\over 16}\right)^{2-\nu_2} 
e^{-i{\pi v(x)\over \omega_1}}
  ~ + O(x^2)
$$
Expanding $v(x)$ and (\ref{uccellogiallo}) we get: 
\be
\label{BBBb}
{
y(x)= x\left[\sin^2 \left(
i{1-\nu_2\over 2}\ln {x\over 16} + {\pi \nu_1\over 2} +\phi(x)\right)  
 ~+\delta^*_E(x)\right] 
,~~~1-\nu_2=\sigma}
\ee
$$
\delta^*_E(x)= \sum_{m_1\geq -1,m_2\geq 1}
a_{m_1m_2}x^{m_1(1-\nu_2)+m_2\nu_2}=O(x),~~~~~
a_{m_1m_2}\in{\bf C}.
$$  
The series converges in ${\cal D}$.


\subsection{ Case $\Re \nu_2=0$, i.e $\Re\sigma=1$}

We observe that $v(x)$ does not vanish when ${\cal V}=0$, because 
$x^{\nu_2}\not\to 0$ . Namely 
$$ 
v(x)= \psi(x)+O(x), ~~~~~
\psi(x):= \sum_{m\geq 1} c_{0m}\left[ e^{i\pi \nu_1} 
x^{\nu_2}\right]^m  \not \to 0  \hbox{ as } x\to 0
$$
The series of $\psi(x)$ converges in ${\cal D}$. We keep the term $\sin^{-2}(f/2)$ and immediately
compute:
$$
  y(x)= \left[{1\over \sin^2\left( -i{\nu_2\over 2} \ln{x\over 16} +{\pi \nu_1\over
  2}+\psi(x)+O(x)\right)}
+O(x^2)\right]
\bigl(1+O(x) \bigr)~+{x\over 2} +O(x^2)
$$
$$
= \left[ \sin^2\left( -i{\nu_2\over 2} \ln{x\over 16} +{\pi \nu_1\over
  2}+\psi(x)\right)+O(x)\right]^{-1}~(1+O(x)),~~~\nu_2=1-\sigma.
$$
If we perform a more explicit computation from (\ref{uccellogiallo})
and the expansion of $v(x)$, 
we get: 
$$
y(x)= \left[ \sin^2\left( -i{\nu_2\over 2} \ln{x\over 16} +{\pi \nu_1\over
  2}+\psi(x)\right)+\sum_{m_1\geq 1} \sum_{m_2\geq -1} A_{m_1m_2}
  x^{m_1(1-\nu_2)+m_2\nu_2}\right]^{-1}
~\times~
$$
$$
\times 
\Bigl(1+ \sum_{m_1\geq 1}\sum_{m_2\geq 0}D_{m_1m_2}
x^{m_1(1-\nu_2)+m_2\nu_2}
\Bigr)
$$
where 
$A_{m_1m_2},B_{m_1m_2},D_{m_1m_2}\in{\bf C}$. 
 The denominator does not vanish on ${\cal D}$. The series are
convergent in ${\cal D}$. 
We can also apply the  the symmetry transformation (\ref{infi}) to
(\ref{BBBb}) and  obtain:
\be
\label{CCCc} 
{
y(x)=\left\{ 
 \sin^2\left( -i{\nu_2\over 2} \ln{x\over 16} +{\pi \nu_1\over
  2}+\psi(x)\right)+\delta^*_E(x)\right\}^{-1}
}
\ee
where 
$$
\delta^*_E(x)= \sum_{m_1\geq 1,m_2\geq -1}
e_{m_1m_2}x^{m_1\sigma+m_2(1-\sigma)}=O(x),~~~e_{m_1m_2}\in{\bf C}
$$
is a convergent series in ${\cal D}$.

\subsection{Representation of solution with oscillatory
  expansions. The bridge between the elliptic representation and the results of this paper}
\label{oscilla}

The identification of  (\ref{mmfido}) with (\ref{BBBb}), and (\ref{XXX}) with
(\ref{CCCc}) is done as follows.  We
rewrite (\ref{BBBb}) and (\ref{CCCc}) in terms of new integration 
 constants $\sigma=1-\nu_2$ and $\phi_E$, instead of  $\nu_2$,
$\nu_1$ (the substitution is obvious). Thus,  (\ref{BBBb}) is:
\be
\label{ondina}
y_E(x)= x\left[\sin^2\left(i{\sigma\over 2}\ln x +\phi_E+
  \sum_{n\geq 1}c_n(\sigma)[e^{-2i\phi_E}x^\sigma]^n\right)+\delta^*_E(x)\right],
\ee
$$
\Re\sigma=0,~~~|x|<\epsilon,~~~|e^{-2i\phi_E}x^\sigma|<\epsilon.
$$
and (\ref{CCCc}) is: 
\be
\label{oondina}
y_E(x)= \left[\sin^2\left(i{1-\sigma\over 2}\ln x 
 +\phi_E+ 
\sum_{n\geq 1} c_n(\sigma) [e^{-2i\phi_E}x^{1-\sigma}]^n
\right)+\delta^*_E(x) \right]^{-1},
\ee
$$
\Re\sigma=1,~~~|x|<\epsilon,~~~|e^{-2i\phi_E}x^{1-\sigma}|<\epsilon.
$$
On the other hand, we have computed the behaviors:
$$
y(x)= x\left[-2iA\sin^2\left(i{\sigma\over 2}\ln x +{\phi\over
    2}-{\pi\over 4}\right)+iA+B+\delta^*(x)\right],~~~\Re\sigma=0.
$$
$$
y(x)= x\left[-2iA\sin^2\left(i{1-\sigma\over 2}\ln x +{\phi\over
    2}-{\pi\over 4}\right)+iA+B+\delta^*(x)\right]^{-1},~~~\Re\sigma=1.
$$
The two  results  must 
coincide, beeing associated to the same monodromy data. . The coincidence is explained by the fact that one can always 
 find an
oscillatory function $f(x)$ such that:
\be
 -2iA\sin^2\left(
{\nu\over 2} \ln x +{\phi\over 2}-{\pi\over 4}
\right)~+iA+B=\sin^2\left(
{\nu\over 2} \ln x +f(x)
\right),~~~\nu\in{\bf R}. 
\label{guarito}
\ee 
If $f(x)$ admits a series expansion (in a suitable domain of
convergence), then  it must have the following form: 
\be
f(x)=\sum_{n\geq 0} f_n x^{-i\nu x}
\label{maldipancia}
\ee
This is exactly the form of the functions in the argument of
$\sin^2(~..~)$ in (\ref{ondina}) and (\ref{oondina}) (just write $\sigma=
-i\nu$ and $\sigma=1+i\nu$ respectively).  This   proves
the convergence of  (\ref{ddelta3}) and (\ref{ddelta2}).

\vskip 0.4 cm 
 The solution $f(x)$ of (\ref{guarito}) is constructed as follows. Let $\psi= {\phi\over 2}-{\pi \over
 4}$. (\ref{guarito}) becomes the equation:
 $$
e^{4if} 
+2
\Bigl[
iAe^{2i\psi}+(2B-1)x^{-i\nu}+iAe^{-2i\psi}x^{-2i\nu}
\Bigr] 
~
e^{2if}
~+x^{-2i\nu}=0
$$
Let $f_1$, $f_2$ be the two solutions: 
 $$
e^{2i(f_1+f_2)}=x^{-2i\nu},
$$
$$
 e^{2if_1}
=
-iAe^{2i\psi}-(2B-1)x^{-i\nu} -iAe^{-2i\psi}x^{-2i\nu}+
$$
$$
-iAe^{2i\psi}
\sqrt{
\left[1+{2B-1\over iA}e^{-2i\psi}x^{-i\nu}+e^{-4i\psi}x^{-2i\nu}\right]^2
+
{1\over A^2}e^{-4i\psi}x^{-2i\nu}
}
$$
The square root is such $-\pi<\arg (\sqrt{...})<\pi$. We observe that 
 $e^{2if_1}$ is clearly an oscillatory function. 
Further observe that the square root is of the form:
$$
\sqrt{
1+ae^{-2i\psi}x^{-i\nu}+be^{-4i\psi}x^{-2i\nu}+ce^{-3i\psi}x^{-3i\nu}+de^{-8i\psi}x^{-4i\nu}
} 
$$ 
where $a,b,c,d$ are constants that can be immediately
computed.
If the absolute value of the sum of the last four terms 
is less then 1 we expand the root in series. 
In particular, this is true if $|e^{-2i\psi}x^{-i\nu}|<r$, for
$r$ suitably small. 
  Thus:
$$
  e^{2if_1}
=
-iAe^{2i\psi}-(2B-1)x^{-i\nu} -iAe^{-2i\psi}x^{-2i\nu}
-iAe^{2i\psi}\Bigl(
1+\sum_{n\geq 1} a_n (e^{-2i\psi} x^{-i\nu})^n 
\Bigr)
$$
$$
 f_1=\psi+{1\over 2i}\ln(-2iA)+{1\over 2i} \ln
\left(
1+{2B-1\over 2iA} e^{-2i\psi}x^{-i\nu}+{1\over 2}
 e^{-4i\psi}x^{-2i\nu}
+{1\over 2}\sum_{n\geq 1} a_n (e^{-2i\psi} x^{-i\nu})^n
\right)
$$
$f_1$ is an oscillatory function. 
If  in a suitable domain  the expansion is possible,
we expand the logarithm and obtain: 
$$
 f_1(x)=\psi +{1\over 2i}\ln(-2iA)~+\sum_{n\geq 1} b_n (e^{-2i\psi} x^{-i\nu})^n
$$
$$
f_2(x)= -f_1-\nu\ln x
$$
Note that the last formula implies: 
$$
\sin^2\Bigl({\nu\over 2} \ln x +f_2\Bigr)= \sin^2\Bigl({\nu\over 2} \ln x +f_1\Bigr).
$$
\qed

\vskip 0.2 cm

\subsection{Example of Picard solutions}  \label{Picard}

 Picard \cite{Picard} studied the case   $\theta_0=\theta_x=\theta_1=0$,
$\theta_\infty=1$. This section is written in order  to show  the general
 results realized in an example that can be computed in terms of
 classical special functions (elliptic and hypergeometric). 
 In this case the function appearing in the elliptic
 representation is $v(x)=0$. Thus: 
$$
  y(x) =\wp\left( \nu_1 \omega_1(x)
  +\nu_2\omega_2(x);\omega_1,\omega^{}_2 \right)+{1+x\over
  3},~~~\nu_1,~\nu_2\in{\bf C},
$$

Apply the Fourier expansion to: 
$$ 
\wp\left( \nu_1 \omega^{}_1(x) +\nu_2\omega^{}_2(x);\omega^{}_1,\omega^{}_2 \right)=\wp\left( \nu_1 \omega^{}_1(x) +[\nu_2+2N]\omega^{}_2(x);\omega^{}_1,\omega^{}_2 \right)
$$
 The domain of convergence is:    
$$ 
\left|\Im \left[{\nu_1\over 2} +\left({\nu_2\over 2} +N\right)\tau(x)  \right] \right|< \Im \tau(x)
$$ 
 Namely: 
\be
   (\Re \nu_2 +2+2N) \ln{|x|\over 16} +O(x)< \Im \nu_2~\arg x +\pi ~\Im \nu_1< 
 (\Re \nu_2 -2+2N) \ln{|x|\over 16}+O(x)
\label{nuovadef}
\ee
This is larger than ${\cal D}$. 
  The critical behavior for $x\to 0$ is computed along the paths:  
$$ 
\arg x =\arg x_0 +{\Re \nu_2 +2N - {\cal V}\over \Im \nu_2} \ln |x|,
~~~ -2 \leq {\cal V} \leq 2,~~~\Im \nu_2\neq 0   
$$
$$
\hbox{If $\Im \nu_2=0$ we take a radial path.}
$$ 

\noindent
 The  critical behavior is then obtained by extracting the leading
 terms of the Fourier expansion. We do this straightforwardly
if $0\leq {\cal V}<2$. The other cases are  obtained from the previous
one by changing $N\mapsto N\pm 1$. Results:  
\vskip 0.2 cm 
\noindent  $\diamond$ 
For $0<{\cal V} <1$ 
 $$
{ 
y(x) = -{1\over 4} \left[ 
{e^{i\pi \nu_1} \over 16^{\nu_2+2N-1}}
\right]
x^{\nu_2+2N} ~(1+O(x^{\nu_2+2N},x^{1-\nu_2-2N}))
}
$$

\vskip 0.2 cm 
\noindent  $\diamond$ 
For $1<{\cal V} <2$ 
 $$
{ 
y(x) = -{1\over 4} \left[ 
{e^{i\pi \nu_1} \over 16^{\nu_2+2N-1}}
\right]^{-1}
x^{2-\nu_2-2N} ~(1+O(x^{2-\nu_2-2N},x^{\nu_2+2N-2}))
}
$$
\vskip 0.2 cm 
\noindent  $\diamond$ 
For ${\cal V}=1$ 
$$
{
y(x) = x~\left[ \sin^2\left(i{1-\nu_2-2N\over 2} \ln{x\over 16} +{\pi \nu_1
\over 2}\right)~+O(x)\right]
}
$$

\vskip 0.2 cm 
\noindent  $\diamond$ 
For ${\cal V}=0$
$$
y(x)=\left[{1
\over  
\sin^{2} \left( 
-i{\nu_2+2N\over 2} \ln{x\over 16} +{\pi \nu_1\over 2} - i {\nu_2+2N\over 2} 
\left[{F_1(x) \over F(x)}+\ln 16\right]
\right)}~+O(x^2)\right]  \times~~~~~~~~~~~~~~~~~~~~~~~~~~
$$
$$
~~~~~~~~~~~~~~~~~~~~~~~~~~~~~~~~~\times \left(1-{x\over 2}+O(x^2)\right)+
~{x\over 2}+O(x^2)
$$
Namely:
$$
{
y(x)=   
\sin^{-2} \left( 
-i{\nu_2+2N\over 2} \ln{x\over 16} +{\pi \nu_1\over 2} - i {\nu_2+2N\over 2} 
\left[{F_1(x) \over F(x)}+\ln 16\right]
\right) ~(1+O(x)) ~+O(x)
}
$$
$$
 =\left[
\sin^{2} \left( 
-i{\nu_2+2N\over 2} \ln{x\over 16} +{\pi \nu_1\over 2}
\right)~+O(x)\right]^{-1}(1+O(x))+O(x)
$$

\vskip 0.2 cm 
\noindent  $\diamond$ 
For ${\cal V}=2$: 
$$
  y(x)=\left[ 
\sin^2 \left( 
i{2-\nu_2-2N\over 2} \ln{x\over 16} +{\pi \nu_1\over 2} 
\right)~+O(x)
\right]^{-1}(1+O(x))+O(x)
$$

\vskip 0.2 cm 
\noindent  $\diamond$ 
For $-1<{\cal V}<0$: behavior of case  $1<{\cal V}<2$ with  $N \mapsto N+1$.

\vskip 0.2 cm 
\noindent  $\diamond$ 
For $-2<{\cal V}<-1$: behavior of case $0<{\cal V}<1$ with  $N \mapsto N+1$.

\vskip 0.2 cm 
\noindent  $\diamond$ 
For ${\cal V}=-1$: behavior of case ${\cal V}=1$ with $N \mapsto N+1$.

\vskip 0.2 cm 
\noindent  $\diamond$ 
For ${\cal V}=-2$: behavior of case  ${\cal V}=0$ with $N \mapsto N+1$.

\vskip 0.2 cm 
\noindent 
$\diamond$ 
If $\Im \nu_2=0$, we choose the convention $ 
 0 \leq \nu_i <2
$. 
The critical behavior for $0<\nu_2<1$ is the same of the case $\Im \nu_2
\neq 0$ with $N=0$ and  $0<{\cal V}<1$; for $1<\nu_2<2$ it 
is the same of the case $\Im \nu_2\neq 0$ with $N=0$ and  $1<{\cal
  V}<2$. Finally, in special cases we have Taylor expansions:
$$
y(x) = x~\left[\sin^2\left({\pi \nu_1\over 2}\right)+\sum_{n\geq 1}a_n
  x^n\right],~~~
\hbox{ if } \nu_2=1
$$
$$
    y(x) = \sin^{-2}\left({\pi \nu_1 \over 2}
\right) +\sum_{n\geq 1} a_n x^n,~~~\hbox{ if }\nu_2=0,~~~\nu_1\neq 0
$$

 \vskip 0.2 cm 
Observe that the choice of $N$ is arbitrary, therefore the {\it same } transcendent has different critical behaviors on different domains (\ref{nuovadef}) 
specified by different values of $N$.

\vskip 0.3 cm 
\noindent
{\bf Remark:} Note that in the cases ${\cal V}=-2,0,2$, the
 denominator 
$\sin^2(...)$ may vanish in the domain (\ref{nuovadef}). 
Therefore, there may be movable poles.   The
 position of the poles can be determined if we keep $F_1(x)/F(x)$ in
 the argument of $\sin^2(~..~)$ and set $\sin^2(~..~)=0$. 

\vskip 0.3 cm

Now let $N=0$ and  $\nu_2=i\nu$, $\nu\in {\bf R}$. Identify $\sigma =
1-\nu_2$. When ${\cal V}=0$, $x \to 0$ along a radial path $\arg x=$
constant. The behavior becomes:
\be
\label{PippoC}
{
 y(x)=
\sin^{-2} \left( 
{\nu\over 2} \ln{x\over 16} +{\pi \nu_1\over 2} +{\nu\over 2} 
\left[{F_1(x) \over F(x)}+\ln 16\right]
\right)~(1+O(x))~+O(x) 
}
\ee

Let $N=0$ and  $\nu_2=1+ i\nu$, $\nu\in {\bf R}$, and  $\sigma =
1-\nu_2$. When ${\cal V}=1$, $x\to 0$ along a radial path, and the
behavior becomes: 
$$
{
 y(x)=~x~\left[
\sin^{2} \left( 
{\nu\over 2} \ln{x\over 16} +{\pi \nu_1\over 2} +{\nu\over 2} 
\right)~+O(x)\right] 
}
$$
From the above computations, we see that the critical behavior of the 
Picard solutions is in accordance with our general results.

\vskip 0.5 cm 
\noindent
{\bf Acknowledgment.} I would like to thank M. Mazzocco for  stimulating discussions and valuable suggestions during the work which brought to  this paper.  I  thank A. Kitaev for valuable comments on the manuscript. I also thank the anonymous referee whose comments improved the paper.


\end{document}